%% file: main.tex
\documentclass[11pt]{article}
\usepackage{jheppub}
\usepackage[dvipsnames]{xcolor}
\usepackage{amsmath,amssymb}
\usepackage{graphicx}
\usepackage{subcaption}
\usepackage{stackengine}
\usepackage{hyperref}
\usepackage{comment}

\makeatletter
\def\@fpheader{\relax}
\makeatother

\title{Colored Jones Polynomials and the Volume Conjecture}

\author[*,\dagger]{Mark Hughes}
\author[\ddagger]{\!\!, Vishnu Jejjala}
\author[\mathsection]{\!\!, P.\ Ramadevi}
\author[\ddagger]{\!\!, Pratik Roy}
\author[\P]{\!\!, Vivek Kumar Singh}

\affiliation[\,*]{Department of Mathematics, Brigham Young University,
275 TMCB, Provo, UT 84602, USA}
\affiliation[\,\dagger]{Max Planck Institute for Mathematics, Bonn, Germany}
\affiliation[\,\ddagger]{Mandelstam Institute for Theoretical Physics, School of Physics, NITheCS, and CoE-MaSS,\\
University of the Witwatersrand, Johannesburg, WITS 2050, South Africa}
\affiliation[\,\mathsection]{Department of Physics, Indian Institute of Technology Bombay, Powai, Mumbai, 400076, India}
\affiliation[\,\P]{Center for Quantum and Topological Systems (CQTS), NYUAD Research Institute,\\
New York University Abu Dhabi, PO Box 129188, Abu Dhabi, UAE}

\emailAdd{hughes@mathematics.byu.edu}
\emailAdd{v.jejjala@wits.ac.za}
\emailAdd{ramadevi@iitb.ac.in}
\emailAdd{pratik.roy@wits.ac.za}
\emailAdd{vks2024@nyu.edu}

\newcommand{\nn}{\nonumber}
\newcommand{\eref}[1]{(\ref{#1})}

\input{custom.tex}

\abstract{
Using the vertex model approach for braid representations, we compute polynomials for spin-$1$ placed on hyperbolic knots up to $15$ crossings.
These polynomials are referred to as $3$-colored Jones polynomials or adjoint Jones polynomials.
Training a subset of the data using a fully connected feedforward neural network, we predict the volume of the knot complement of hyperbolic knots from the adjoint Jones polynomial or its evaluations with $99.34\%$ accuracy.
A function of the adjoint Jones polynomial evaluated at the phase $q=e^{ 8 \pi i  \over 15 }$ predicts the volume with nearly the same accuracy as the neural network.
From an analysis of $2$-colored and $3$-colored Jones polynomials, we conjecture the best phase for $n$-colored Jones polynomials, and use this hypothesis to motivate an improved statement of the volume conjecture.
This is tested for knots for which closed form expressions for the $n$-colored Jones polynomial are known, and we show improved convergence to the volume.

}

\begin{document}

\maketitle

\section{Introduction and summary}
Chern--Simons theory is arguably the simplest non-trivial quantum field theory.
Natural operators in this three-dimensional topological field theory are Wilson loops:
\begin{equation}
U_R(\Gamma) = \text{tr}_R\; \mathcal{P} \exp \big( i \oint_\Gamma A \big) \,,
\end{equation}
traces of path ordered exponentials of the holonomy of the gauge field along paths $\Gamma$.
The connection $A = A_\mu\, dx^\mu$ is a Lie algebra valued $1$-form and $R$ denotes a representation of the associated gauge group $G$.
Suppose we let $\Gamma$ be a knot $K$, which is an embedding of $S^1$ into a three-manifold, which we take to be $S^3$.
Colored Jones polynomials $J_{n}(K; q)$ are given by the expectation value of the Wilson loop operator along the knot $K$, with $R$ being a $n$-dimensional representation of $SU(2)$~\cite{Witten:1988hf}.
More precisely, we normalize as follows:
\begin{equation}
J_{n}(K; q) := \frac{\int_\mathcal{U} U_n(K)\; e^{iS_\text{CS}}}{\int_\mathcal{U} U_n(\bigcirc)\; e^{iS_\text{CS}}} = \frac{\langle U_n(K) \rangle}{\langle U_n(\bigcirc) \rangle} \,,
\end{equation}
where the Chern--Simons action
\begin{equation}\label{eq:SCS}
S_\text{CS} = \frac{k}{4\pi} \int_{S^3} \big( A\wedge dA + \frac23 A\wedge A\wedge A \big) \,,
\end{equation}
$\mathcal{U}$ is the space of $SU(2)$ connections modulo gauge transformations, and $\bigcirc$ denotes the unknot.
The gauge invariance of~\eqref{eq:SCS} requires the integer quantization of the Chern--Simons level, $k\in \mathbb{Z}$.
The $J_{n}(K; q)$ are Laurent polynomials in a variable $q$, which is related to the Chern--Simons coupling:
\begin{equation}
q := \exp \big( {2 \pi i\over k+2} \big) \,.
\end{equation}
These polynomials are topological invariants of the knot, meaning they are independent of how the knot is drawn and the metric on the underlying three-manifold.
The $n=2$ case was initially obtained from finite dimensional von Neumann algebras~\cite{jones1985polynomial,jones1987hecke} and then interpreted combinatorially~\cite{kauffman1987state}.
Skein relations make manifest that the coefficients in the polynomial are integer valued.
This is also established in Khovanov homology~\cite{khovanov2000categorification,bar2002khovanov}.

A hyperbolic knot is one for which the knot complement, obtained from excising a tubular neighborhood around the knot in $S^3$, admits a complete Riemannian metric of constant negative curvature whose uniqueness is ensured by Mostow--Prasad rigidity.
The volume of the knot complement, $V(S^3\backslash K)$, computed using this metric is a knot invariant~\cite{thurston_geometry_topology_three_manifolds}.
At low crossing number, nearly every knot is hyperbolic, but the proportion of hyperbolic knots does not approach unity in the large crossing number limit~\cite{belousov2019hyperbolicknotsgeneric}.
Thus, in this work and elsewhere, caution is to be applied in extrapolating phenomena from datasets that are inherently atypical.

Kashaev conjectured~\cite{kashaev1995link} that the colored Jones polynomial in the large color ($n\to\infty$) limit is related to $V(S^3\backslash K)$:
\begin{equation}\label{eq:vc}
	\lim_{n\to\I} \frac{\log|J_n(K;\omega_n)|}n = \frac1{2\pi} V(S^3\backslash K) \,, \quad \text{where} \quad \omega_n := e^{2\pi i/n} \,.
\end{equation}
(See also~\cite{murakami2001colored, Gukov:2003na}.)
This volume conjecture relates a quantum invariant, the colored Jones polynomial, to a classical geometric quantity, the volume.
The semiclassical limit where the volume conjecture is realized is, in fact, a double scaling limit~\cite{witten2011analytic}:\footnote{
As explained in~\cite{Craven:2020bdz}, it is both convenient and natural for us to define $\gamma$ in this manner.
The semiclassical limit is alternatively phrased as $q^n = e^{2\pi i \widetilde\gamma}$, which puts $\widetilde\gamma := n/(k+2) = \text{constant}$ as $n\to\infty$~\cite{witten2011analytic}.
}
\begin{equation}
\gamma := \frac{n-1}{k} \to 1 \quad \text{as} \quad n\,,\; k \to\infty \,.
\end{equation}
The analytic continuation of Chern--Simons theory~\cite{witten2011analytic} suggests the possibility of a non-linear relation between the colored Jones polynomial and the hyperbolic volume for every color $n$.\footnote{
See also~\cite{Garoufalidis_2011}, which, exploiting the cyclotomic expansion~\cite{Habiro_2007}, discusses the asymptotics of colored Jones polynomials.}
In the simplest case, such a relation was deduced from a dataset of knots up to $16$ crossings for which fundamental representation Jones polynomials $J_2(K;q)$ have been tabulated~\cite{KnotAtlas,linkinfo,SnapPy}.
Ref.~\cite{Jejjala:2019kio} showed that a deep neural network could predict the hyperbolic volume from the coefficients and the maximum and minimum degrees of the Jones polynomials to better than $98\%$ accuracy.
Since two knots can have the same fundamental representation Jones polynomial and different hyperbolic volumes, and when this happens the volumes differ by about $3\%$, the neural network performs about as well as can be expected.
In~\cite{Craven:2020bdz}, the authors used feedforward neural networks to determine that a function of $|J_2(K; e^{3\pi i/4})|$ approximates the volume with nearly the same accuracy.
Translating $e^{3\pi i/4}$ to a Chern--Simons coupling, we associate the phase in this evaluation to the fractional level $k=\frac23$.

The explanation for why such a formula works so well follows from Witten's analytic continuation of Chern--Simons theory~\cite{witten2011analytic}.
In the saddle point approximation to the path integral of the analytically continued Chern--Simons theory, there are particular flat $\mathfrak{sl}(2,\mathbb{C})$ connections of interest.
Heuristically, such geometric connections yield a semiclassical contribution of the form
\begin{equation}
Z \sim e^{iS(\mathcal{A}_+)} \big( 1 - e^{2\pi i k} \big)
\end{equation}
to the partition function.
The approximation formula applies for levels $k$ for which $\mathcal{A}_+$ makes a contribution to the Chern--Simons path integral for a large fraction of knots~\cite{Craven:2020bdz}.
Phases corresponding to integer level do not have access to this information, so the analytic continuation is crucial to recovering the volume.

As the dimension of the $SU(2)$ representation $n$ becomes large, the higher colored Jones polynomials are expected to discriminate knots with different hyperbolic volumes.
\textit{A priori}, especially for small $n$, it is not clear whether this improvement in performance occurs each time we increment the color.
Closed form expressions or formul\ae\ in terms of $q$-Pochhammer symbols exist for colored Jones polynomials for some families of knots like torus knots $T(2,2p+1)$ and twist knots $K_p$,
for instance.  The well-known lowest crossing hyperbolic knot, called the figure-eight knot, $4_1\equiv K_{-1}$, has colored Jones polynomial  given in~\cite{Le:2000, Habiro:2000}.
Another knot for which we can readily compute colored Jones polynomials for any $n$ is $K_0$, which is defined as the closure of the $3$-strand braid $\sigma_1^2 (\sigma_1 \sigma_2)^8$ and was studied in~\cite{Garoufalidis:2004}.
This is a knot with $18$ crossings, but a tetrahedral decomposition of its complement is accomplished with only four ideal tetrahedra.
In this case, numerics show that~\eref{eq:vc} does indeed converge to the volume, but this convergence is slow and non-monotonic. It is so far an open question whether this behavior is characteristic or an aberration.

There are $313{,} 210$ hyperbolic knots up to fifteen crossings~\cite{hoste1998first}.
In this work, we construct a partial dataset of $177 {,}316$ adjoint representation ($n=3$) Jones polynomials for these knots using the vertex model~\cite{Akutsu:1987dz,Ramadevi2017}.
These correspond to the closures of $m$-stranded braids for $m\le 7$.
The error in the neural network prediction of the volume of the knot complement using the adjoint polynomial invariants drops to $0.4\%$.
The approximation formula here uses an evaluation of the $n=3$ colored Jones polynomial at the phase $e^{8\pi i/15}$, corresponding to Chern--Simons level $k=\frac74$.
Based on these experiments, we conjecture the formul\ae
\begin{equation}        \label{eq:best-phase-guess}
q(n) = \exp\left( 2\pi i \frac{n+1}{n(n+2)} \right) \,, \qquad k(n) = \frac{n^2-2}{n+1}
\end{equation}
as the relevant phase and fractional level as a function of $SU(2)$ representation.

This paper employs machine learning to correlate the new dataset of $3$-colored Jones polynomials to the volume of the knot complements of hyperbolic knots.
After~\cite{hughes2016neural}, the utility of machine learning as a tool in low dimensional topology is by now well established.
Representative works in this direction include~\cite{Jejjala:2019kio,levitt2022big,Gukov:2020qaj,Craven:2020bdz,dlotko2023mappertypealgorithmscomplexdata,davies2021advancing,Craven:2021ckk,Craven:2022cxe,Gukov:2023kvx,Gukov:2024buj,Gukov:2024opc}.

The organization of this paper is as follows.
In Section~\ref{sec:data}, we discuss the dataset of $3$-colored Jones polynomials that we generate.
The method for constructing these is reviewed in Appendix~\ref{sec:app}.
In Section~\ref{sec:analysis}, we examine features of $3$-colored Jones polynomials and compare with $2$-colored Jones polynomials.
In particular, we look at their zeros and the statistics of the degrees, the lengths of the polynomials, and evaluations at special points motivated by the volume conjecture.
In Section~\ref{sec:ml}, we apply machine learning to the dataset of $3$-colored polynomials using the degrees and coefficients and using evaluations at phases as inputs to a neural network.
In Section~\ref{sec:newvc}, from data for Jones polynomials in the fundamental and adjoint representations of $SU(2)$, we make a hypothesis for the best phase for an $n$-dimensional representation of $SU(2)$ and test an improved statement of the volume conjecture with known formul\ae\ for the colored Jones polynomials for $4_1$ and $K_0$.
Our code is available on \href{https://github.com/roypratik92/ColoredJonesML}{GitHub}~\cite{github}, and the data are available on \href{https://doi.org/10.5281/zenodo.1490093}{Zenodo}~\cite{data}.

\section{Comments about data}\label{sec:data}

There are a total of $12{,} 965$ knots up to $13$ crossings, of which ten are non-hyperbolic.
At $14$ crossings, there are $46{,} 969$ hyperbolic knots, at $15$ crossings there are $253{,} 285$ hyperbolic knots, and at $16$ crossings there are $1{,} 388{,} 694$ hyperbolic knots.
In total, we have $1{,} 701{,} 903$ hyperbolic knots up to $16$ crossings~\cite{hoste1998first}.
Our dataset has $J_2$ for all of these.
Several knots, even those whose complements have different volumes, can have the same Jones polynomial.
There are $841{,} 139$ unique polynomials in the dataset.

For $J_3$, our dataset only contains the polynomials of knots which are the closures of $m$-strand braids, where $m\leq 7$.
Therefore, up to $13$ crossings, we have $J_3$ for $11{,} 941$ knots.
At $14$ crossings, we have $J_3$ for $18{,} 353$ hyperbolic knots.
At $15$ crossings, we have $J_3$ for $147{,} 022$ hyperbolic knots.
The results presented below thus use a dataset with a total of $177,316$ $J_3$ polynomials for hyperbolic knots up to $15$ crossings of which $164,455$ are unique.
It is difficult to extrapolate based on an incomplete dataset, but a larger fraction of knots are distinguished by $J_3$ in comparison to $J_2$ for the adjoint polynomials we can calculate.
This observation aligns with the expectation of the volume conjecture that larger colors should split degeneracies in the lower colored data.
This information is summarized below in Table~\ref{table:knotdata}.

\begin{table}[h!]
\centering
\begin{tabular}{c|cccc}
\hline
\textbf{Crossings} & \textbf{Total knots} & \textbf{Hyperbolic knots} & \textbf{$J_2$ Computed} & \textbf{$J_3$ Computed} \\
\hline
$\leq 13$ & 12{,}965 & 12{,}955 & 12{,}955 & 11{,}941 \\
14        & 46{,}972      & 46{,}969  & 46{,}969  & 18{,}353 \\
15        & 253{,}293      & 253{,}285 & 253{,}285 & 147{,}022 \\
16        & 1{,}388{,}705      & 1{,}388{,}694 & 1{,}388{,}694 & 0 \\
\hline
\end{tabular}
\caption{Summary of total knots, hyperbolic knots, and the numbers of knots for which $J_2$ and $J_3$ are known, organized by crossing number.}
\label{table:knotdata}
\end{table}

For our machine learning experiments, we use the $1{,} 701{,} 903$ $J_2$ polynomials and the $177{,} 316$ $J_3$ polynomials when considering them separately.
When we need to compare $J_2$ and $J_3$, or when we use their combination for machine learning, we only use the knots for which we know both $J_2$ and $J_3$.  

Note that since we are mostly interested in phases between $0$ and $e^{i\pi}$, we express the phases as $e^{2\pi ix}$ using the variable $x\in[0,\frac12]$ in what follows, including in the plots.

\begin{figure}[h]
\centering
\includegraphics[width=.8\textwidth]{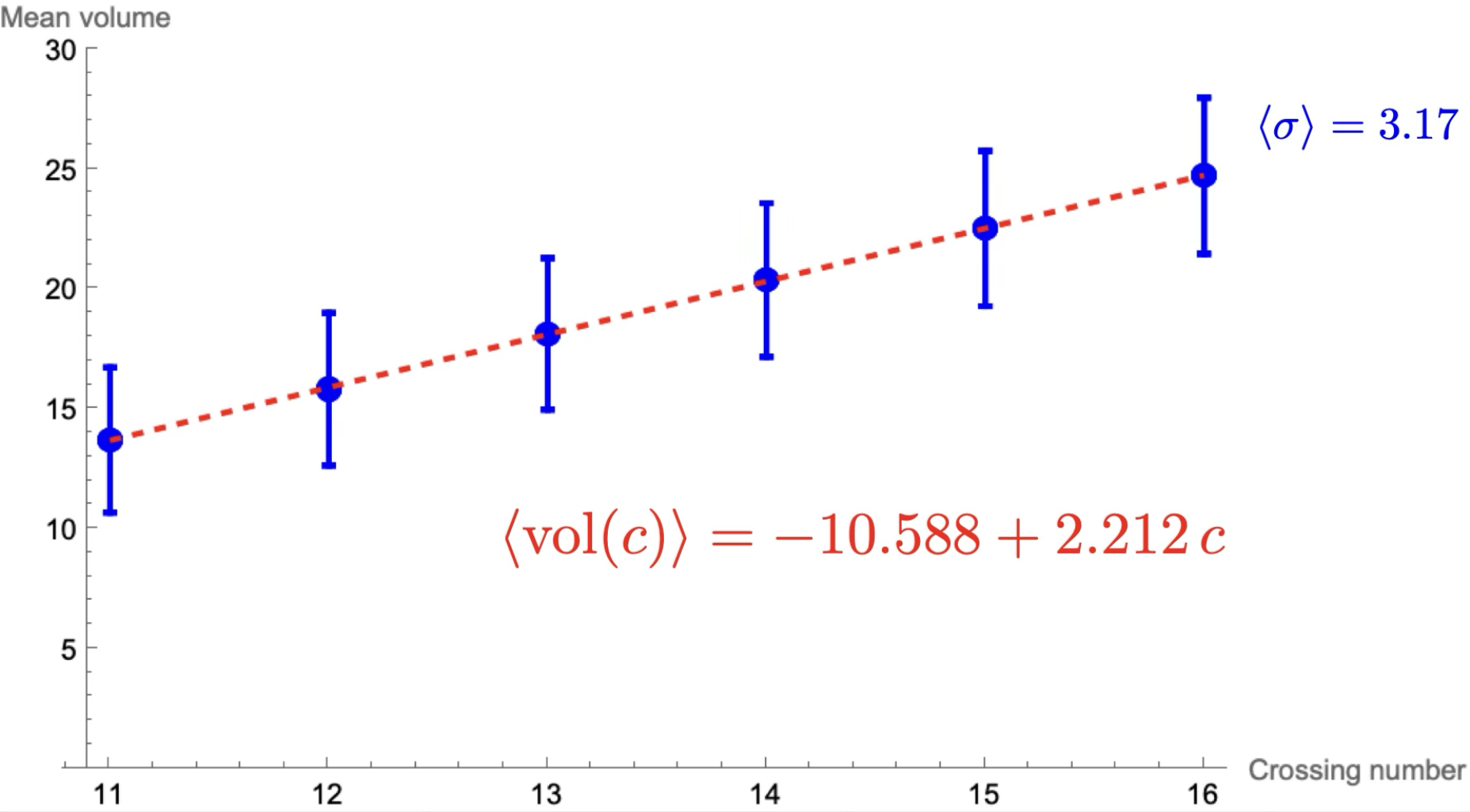}
\caption{\textsf{Mean volume as a function of crossing number.}} 
\label{fig:volumes}
\end{figure}

From the beginning, knot tables have listed knots by their crossing number $c$, and certainly crossing number is a reasonable proxy for the complexity of a knot.
It is far from clear, however, whether this is the best organizing principle there is.
In the case of the volume, it appears to be a useful one.
Figure~\ref{fig:volumes} shows a linear growth of the volume of the knot complement of a hyperbolic knot with the crossing number.
The variance is approximately constant independent of crossing number.
A hint for why this is the case may be found in~\cite{yokota2000volume}, where volumes are described in terms of ideal tetrahedra at the crossings.\footnote{
We thank Sergei Gukov for a discussion about this and for pointing us to the relevant literature.}
%\VJ{Is this the correct reference?}

\section{Numerical analysis of the polynomial data}\label{sec:analysis}

\subsection{Zeros of the polynomials}
The zeros of polynomials, such as the Alexander and Jones polynomials, provide deep insights into the structure of knots and links.
In 2022, Hoste, based on computer experiments, conjectured that for alternating knots, the real part of any zero of the Alexander polynomial satisfies $\operatorname{Re}(z) > -1$.
The conjecture was proven for certain knot classes~\cite{Ishikawa2019, Stoimenow2011,alsukaiti2024alexander,HM2013}, though counterexamples to the general conjecture were found~\cite{HIRASAWA201948}, prompting an exploration of the distributions of zeros of Alexander polynomials.
Similarly, the zeros of Jones polynomials reveal interesting patterns.
For instance, they are dense in the complex plane for certain links, such as pretzel links~\cite{Jin2010}, though they tend to cluster around the unit circle for other classes of links~\cite{AI2005,XF2010}.
Motivated by these results and others~\cite{Andersen:2004}, we study the distribution on the complex plane of zeros of Jones polynomials.

\begin{figure}[h] 
  \begin{center}
    \includegraphics[width=\textwidth]{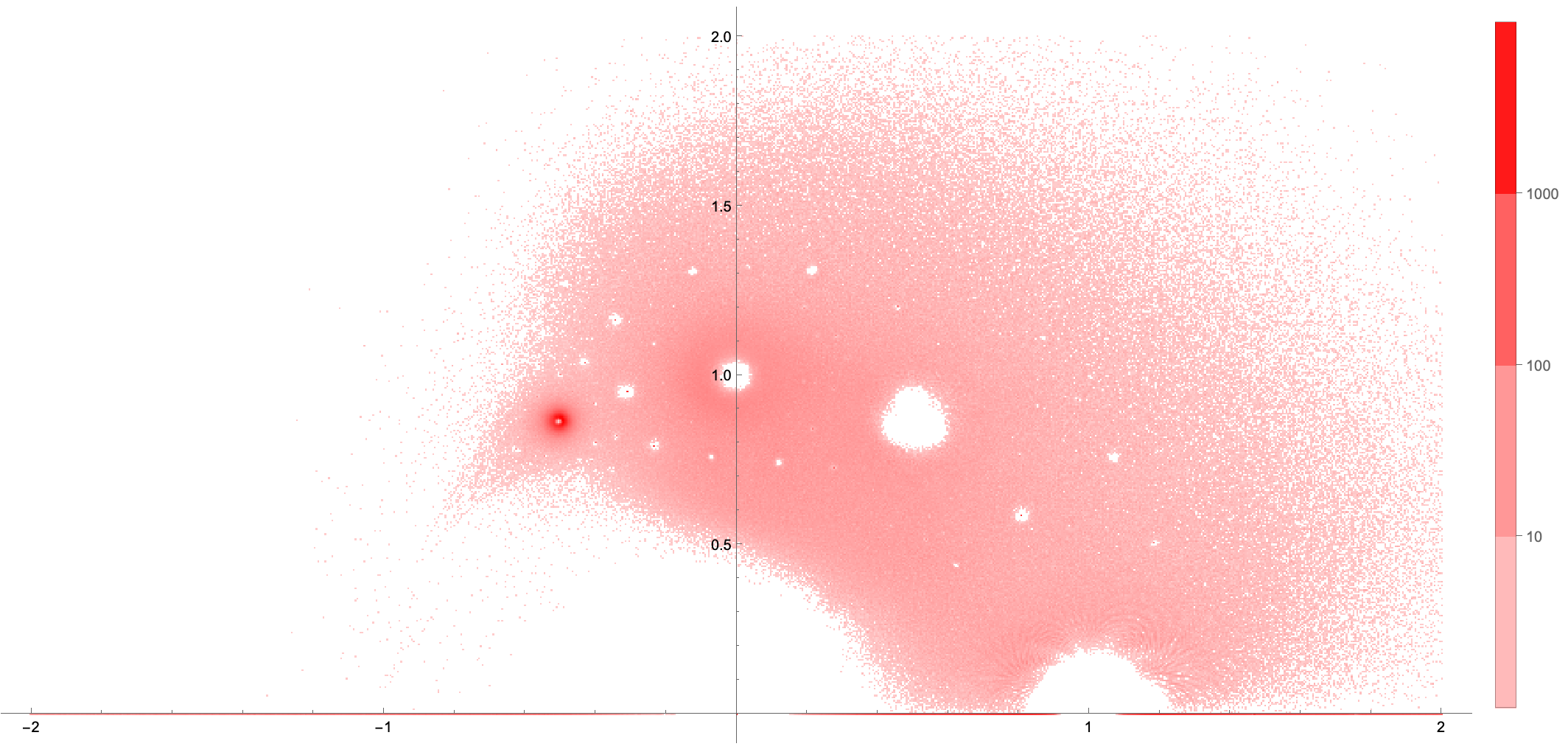}
  \end{center}
  \caption{\textsf{Zeros of $J_2$ polynomials in our dataset in the (upper-half) complex $q$-plane.}}  \label{fig:J2zeros}
\end{figure}

In Figure~\ref{fig:J2zeros}, we plot the zeros of the $1{,}701{,}903$ fundamental Jones polynomials in our dataset on the complex plane.
We have restricted the range of the plot so that the real parts of the roots lie in the interval $(-2,2)$ and the imaginary parts lie in the interval $[0,2)$.
However, there aren't many zeros outside of the plotted range except those that lie on the real axis. 

\begin{figure}[h] 
  \begin{center}
    \stackinset{l}{4pt}{t}{4pt}{\frame{\includegraphics[width=.35\textwidth]{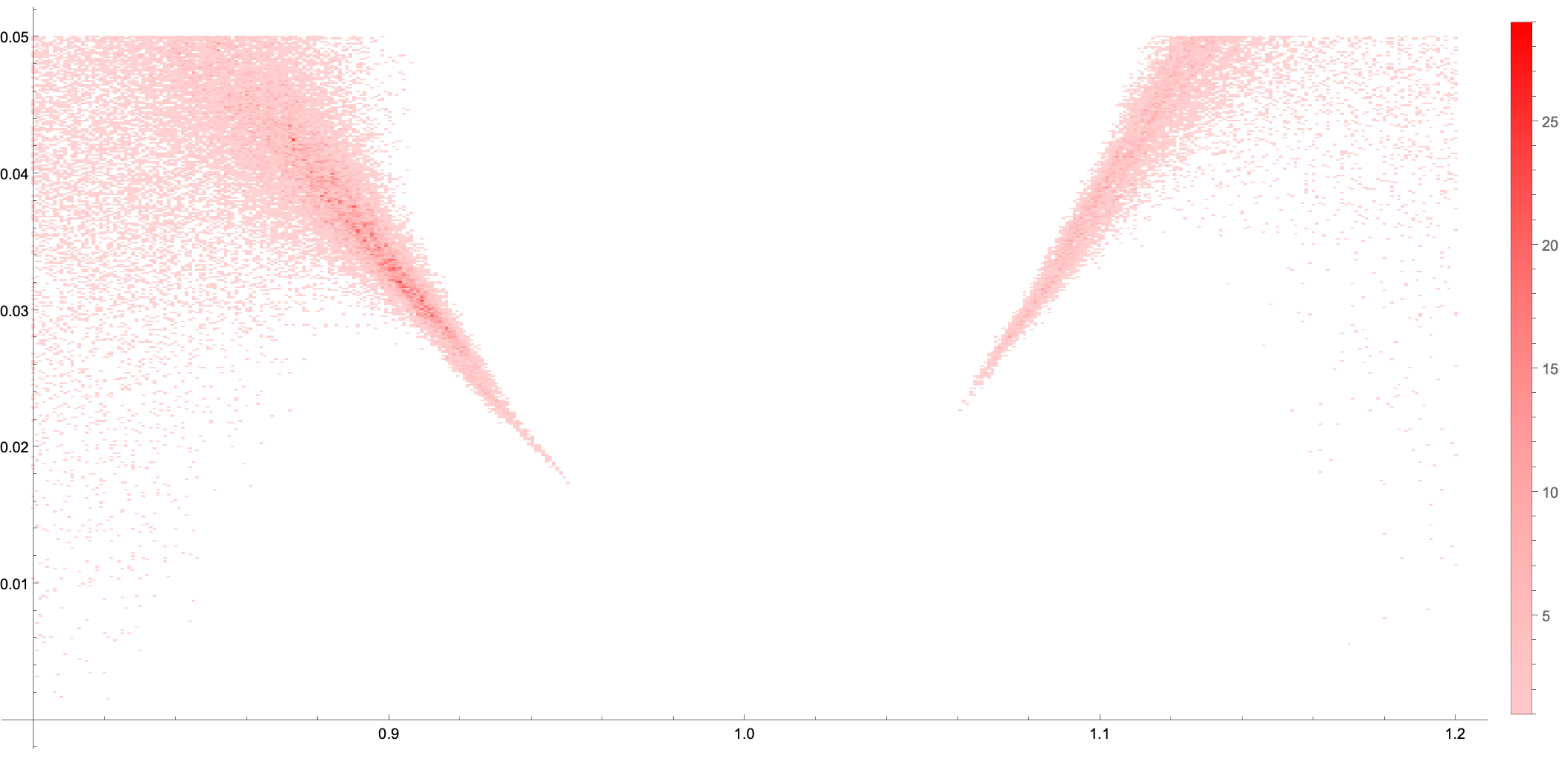}}}{\includegraphics[width=\textwidth]{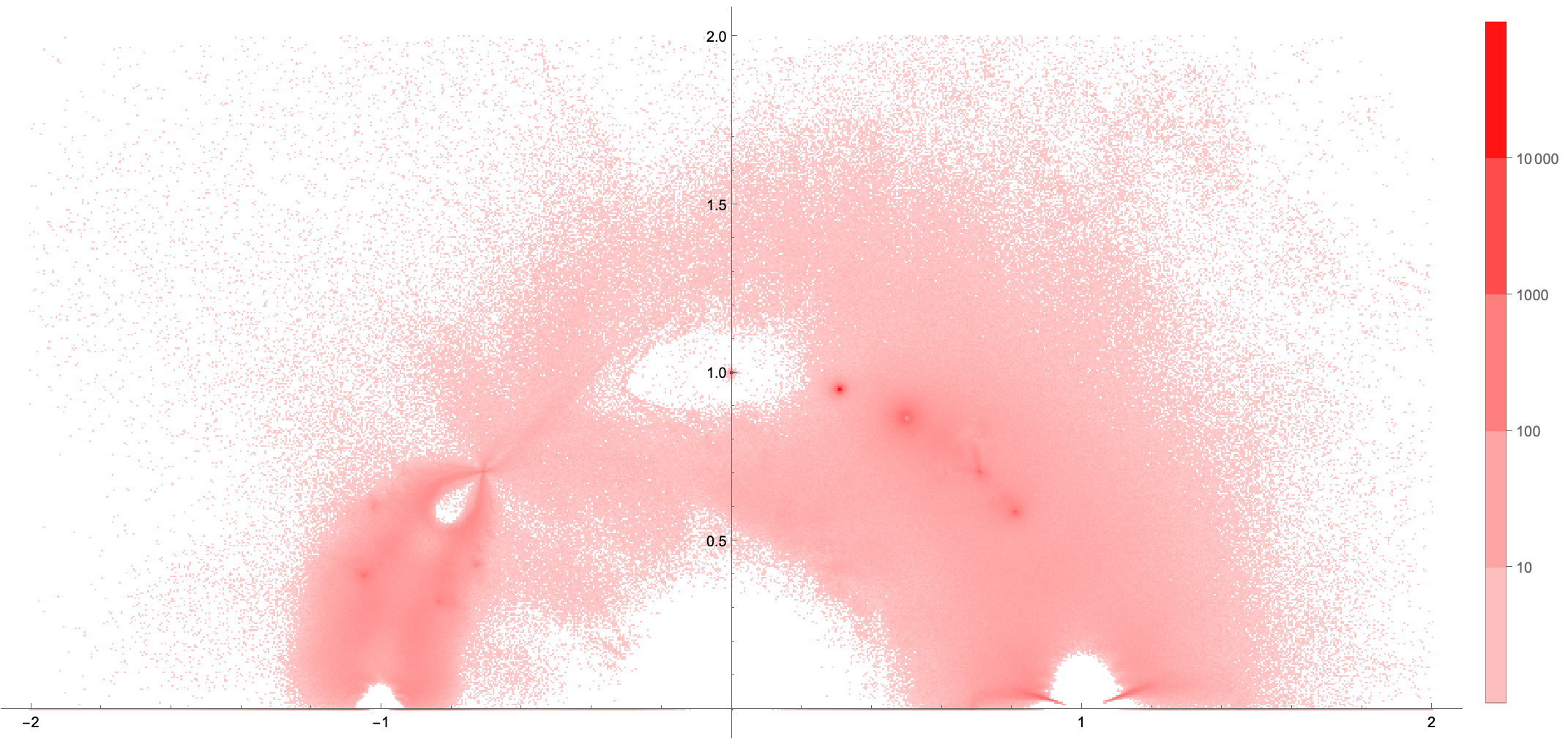}}
  \end{center}
  \caption{\textsf{Zeros of $J_3$ polynomials in our dataset in the (upper-half) complex $q$-plane. Inset: Zoom in of the plot to the $x$-axis interval $(0.8,1.2)$ and $y$-axis interval $(0.0,0.05)$.}}  \label{fig:J3zeros}
\end{figure}

In Figure~\ref{fig:J3zeros}, we plot the zeros of the $177{,} 316$ adjoint Jones polynomials in our dataset on the complex plane.
The range of $z$ in Figure~\ref{fig:J3zeros} is restricted as in Figure~\ref{fig:J2zeros}.
Note that in both Figures~\ref{fig:J2zeros} and~\ref{fig:J3zeros}, there are a lot of zeros at $0$; these are not visible due to our methods of plotting.
In particular, $28{,} 345$ of the fundamental Jones polynomials of the $313{,} 209$ hyperbolic knots up to $15$ crossings have zeros at zero.
Similarly, $36{,} 918$ of the adjoint Jones polynomials of $177{,} 316$ hyperbolic knots in the dataset have zeros at zero.

In the inset at top left of Figure~\ref{fig:J3zeros}, we zoom in on one feature of the zeros of $J_3$.
We find that the $J_3$ polynomials in our dataset do not have any zeros in a neighborhood of $1$, but there is some interesting substructure there.
(The Jones polynomial never has a zero at $1$.)

%\PR{Maybe add comments about the distributions of zeros of polynomials in general, in comparison with, e.g.,~\cite{Mezincescu:1996, Tao:2013}.}


\subsection{Degrees of polynomials}
\begin{figure}[h]
\begin{tabular}{cc}
  \subcaptionbox{\textsf{Minimum degrees of $J_2$ and $J_3$.}}{\includegraphics[width=75mm]{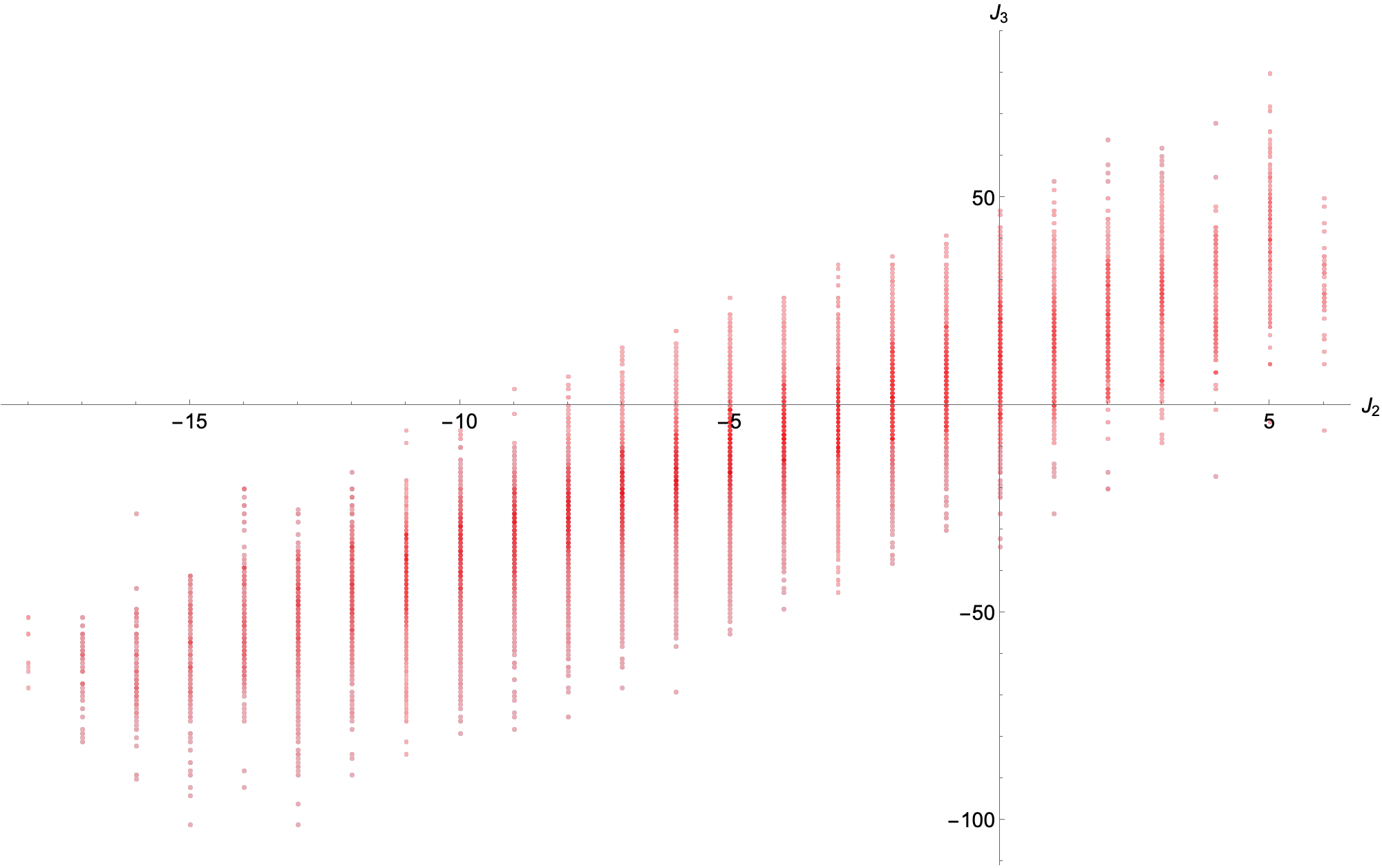}}  &   
  \subcaptionbox{\textsf{Maximum degrees of $J_2$ and $J_3$.}}{\includegraphics[width=75mm]{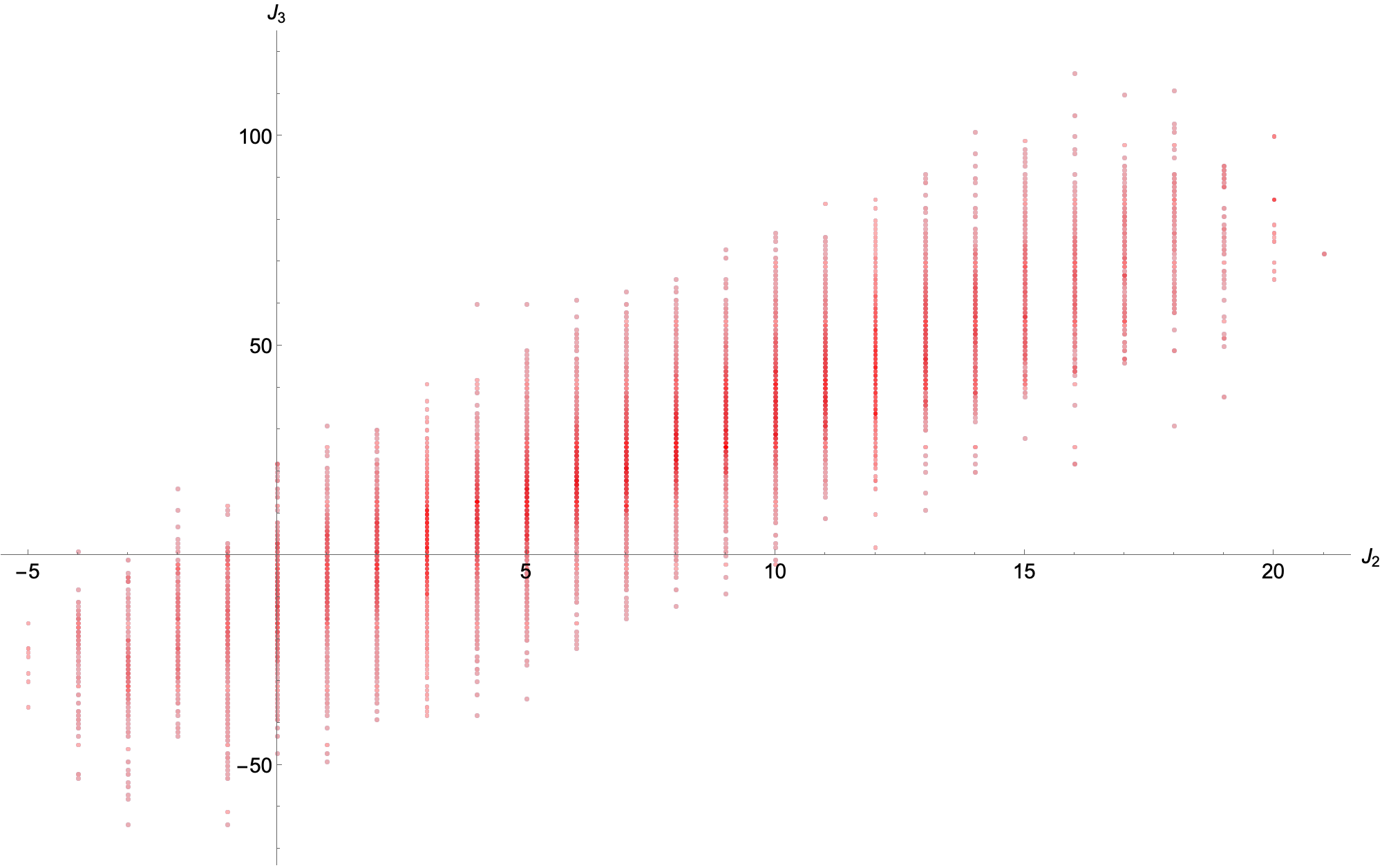}} \\[15pt]
  \subcaptionbox{\textsf{Minimum and maximum degrees of $J_2$.}}{\includegraphics[width=75mm]{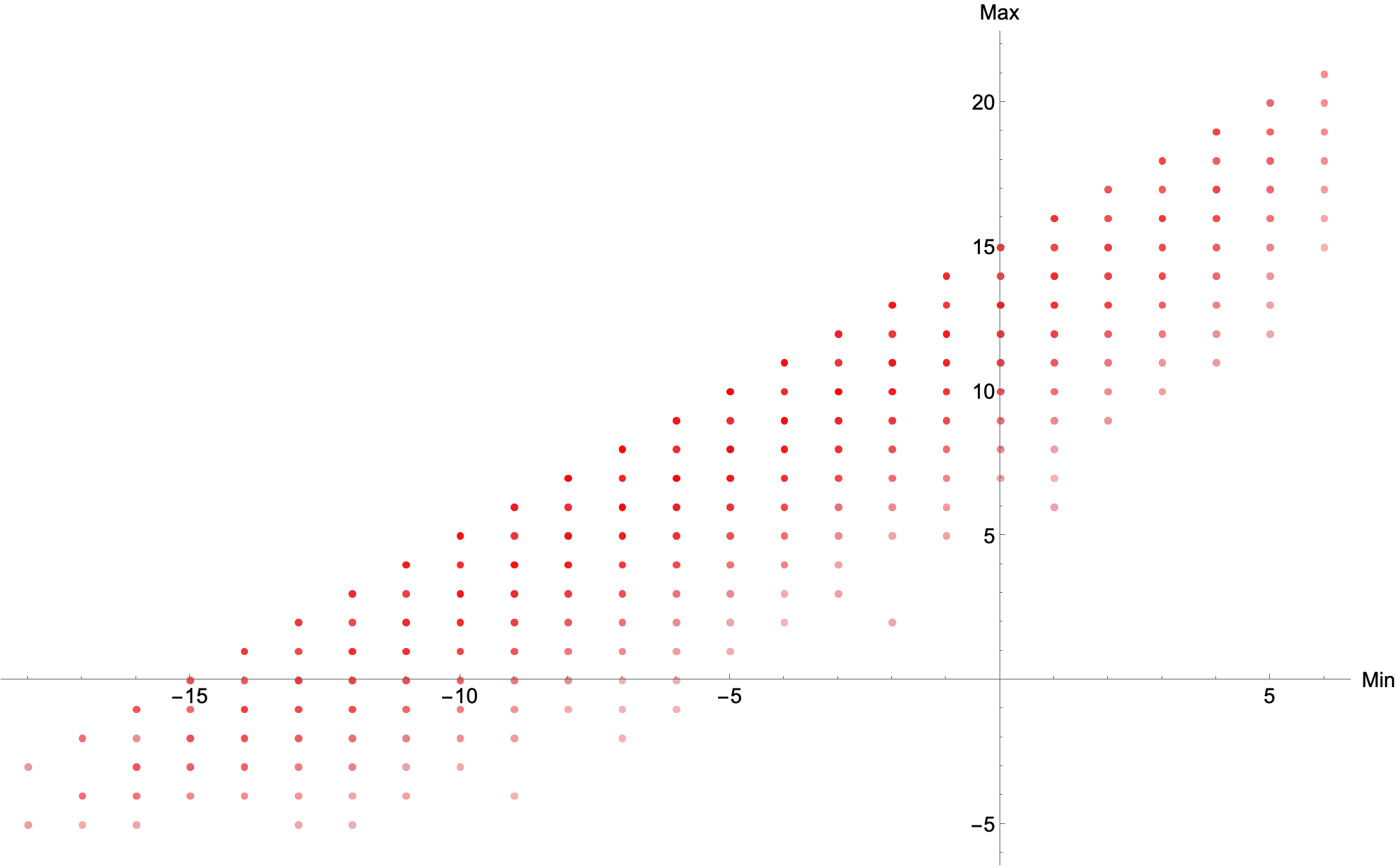}} &    \subcaptionbox{\textsf{Minimum and maximum degrees of $J_3$.}}{\includegraphics[width=75mm]{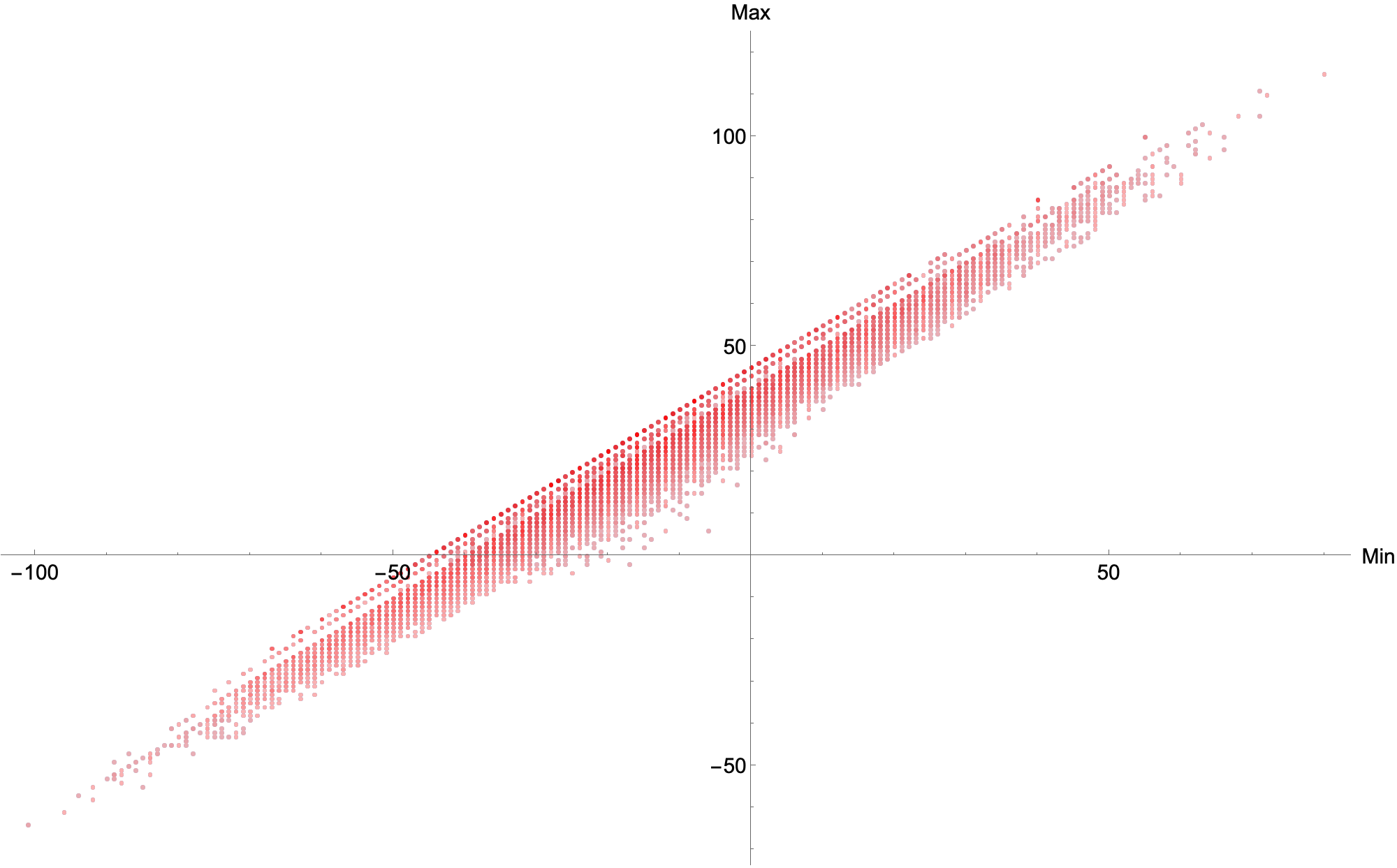}} 
\end{tabular}
\caption{\textsf{Correlations between degrees of polynomials. Darker shades on points mean higher frequency.}}
\label{fig:j2j3minmaxdegs}
\end{figure}

We also study the distribution of the maximum and minimum degrees of the Jones polynomials.
Figure~\ref{fig:j2j3minmaxdegs} shows the distributions of the minimum and maximum degrees of the polynomials and the correlations between minimum and maximum degrees for the different colors.

\begin{figure}[h] 
	\centering
	\begin{minipage}{0.45\textwidth}
  	  \begin{center}
    	  \includegraphics[width=\textwidth]{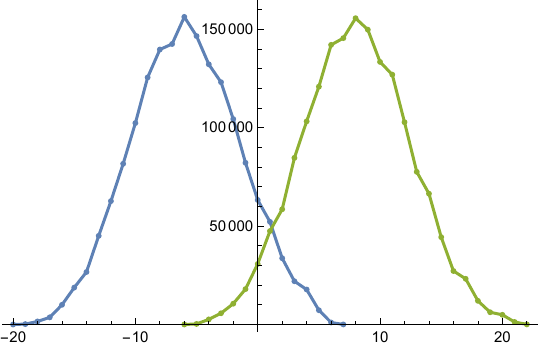}
  	  \end{center}
	\end{minipage}
	\begin{minipage}{0.45\textwidth}
  	  \begin{center}
    	  \includegraphics[width=\textwidth]{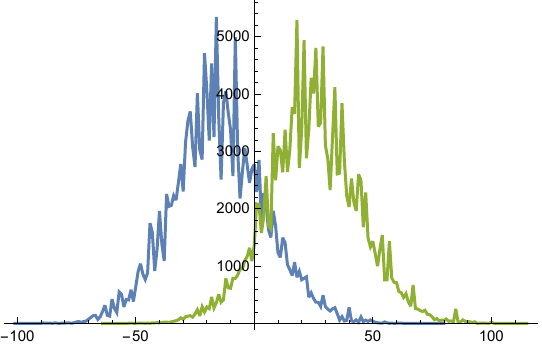}
  	  \end{center}
	\end{minipage}
  \caption{\textsf{Distribution of minimum (blue) and maximum (green) degrees of $J_2$ (left subfigure) and $J_3$ (right subfigure). The $y$-axis corresponds to the number of times the given degree occurs as the minimum or maximum in our dataset.}}  \label{fig:j2j3minmax-dist}
\end{figure}

We find that the minimum and maximum degrees of the polynomials have a roughly Gaussian distribution, as seen in Figure~\ref{fig:j2j3minmax-dist}.

\begin{figure}[h] 
  \begin{center}
    \includegraphics[width=.7\textwidth]{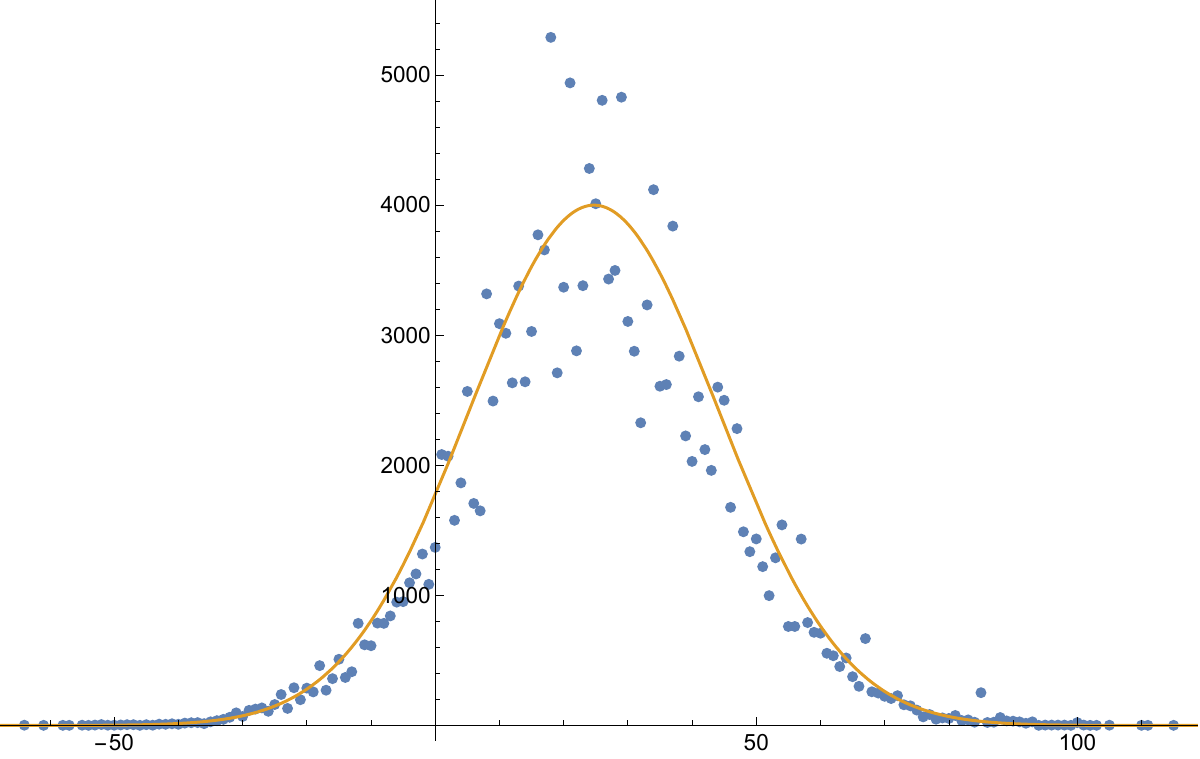}
  \end{center}
  \caption{\textsf{Plot of maximum degrees of $J_3$ with a Gaussian fit in yellow for comparison.}}  \label{fig:j3max-fit}
\end{figure}

In Figure~\ref{fig:j3max-fit}, we plot a Gaussian function with mean and variance determined by the distribution of the maximum degrees of $J_3$ polynomials as an example.

We also note that the minimum and maximum degrees of the polynomials grow at a roughly similar rate, at least within our dataset.
The means of the minimum and maximum degrees of $J_2$ are $-5.86$ and $8.04$, and the modes are $-6$ and $8$.
For $J_3$ the minimum and maximum degrees have means $-15.03$ and $24.73$, and the modes are $-16$ and $18$.
The ratio of the mean minimum degrees of $J_3$ and $J_2$ is $2.56$.
The ratio of the mean maximum degrees of $J_3$ and $J_2$ is $3.07$.

\begin{figure}[h]
\centering
\includegraphics[width=.8\textwidth]{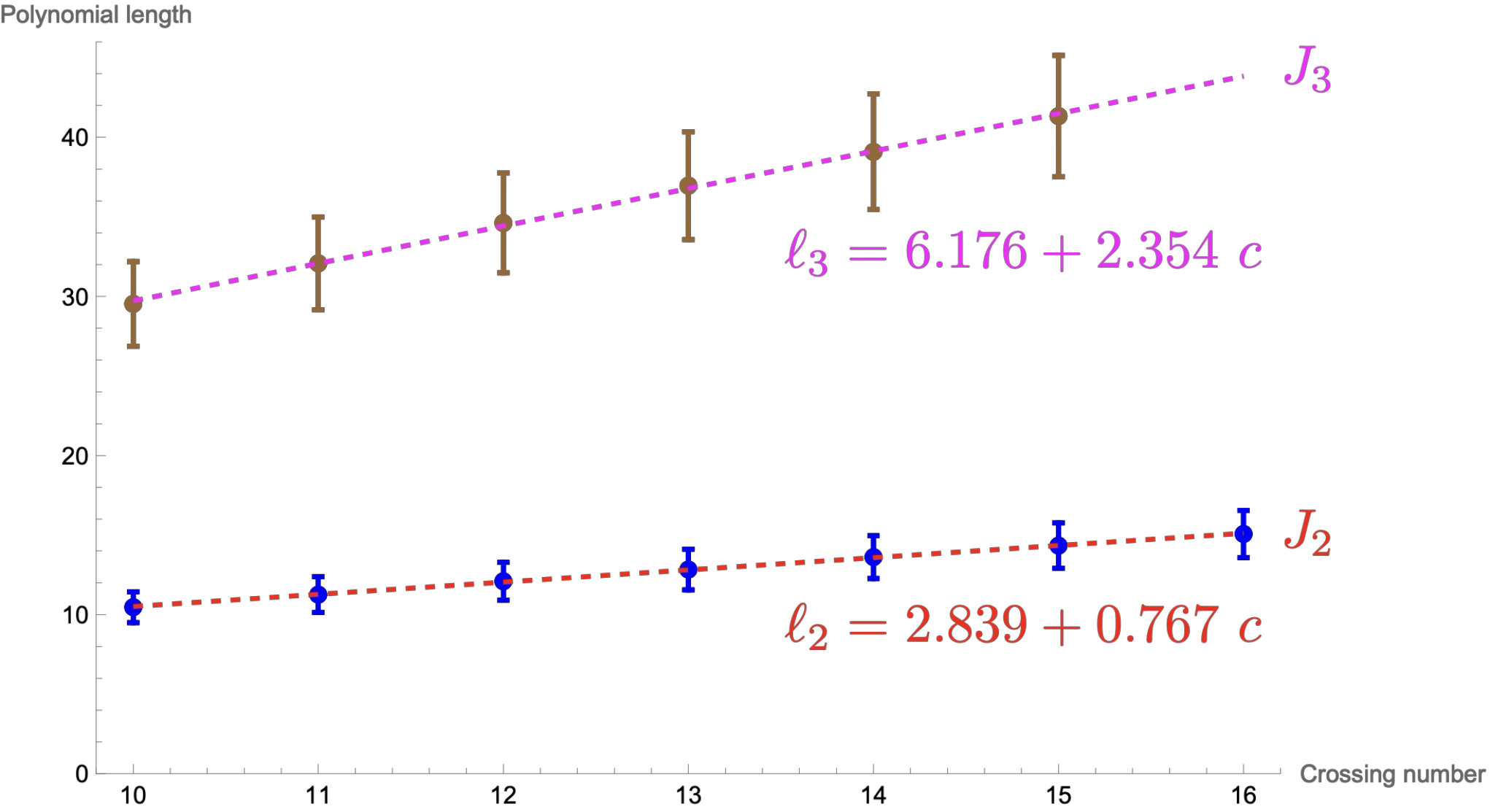}
\caption{\textsf{Length of fundamental and colored Jones polynomials as a function of crossing number.}} 
\label{fig:length}
\end{figure}

As seen in Figure~\ref{fig:length}, the length of the Jones polynomial (\textit{i.e.}, the maximum degree minus the minimum degree plus one) follows a linear relationship with crossing number.
The standard deviations are computed using the set of available polynomials for the knots in the dataset at each crossing number.
The linear fit may be motivated by the calculation of the Jones polynomials in terms of traces of $R$-matrices defined at each of the crossings.
Normalizing appropriately, this gives a standard length.
To deviate from this standard length requires, \textit{e.g.}, detailed cancellations, and this is generically rare.\footnote{
We thank Sergei Gukov for discussing this result with us and suggesting this explanation.}
If one were to assign to every knot the mean volume of all the knots in the dataset, the error in the prediction is about $12\%$.
The linear relationships in Figures~\ref{fig:volumes} and~\ref{fig:length} and the correlations therein may describe the next correction, and machine learning the volume conjecture is a refinement of this.

\subsection{Evaluation of polynomials}

\begin{figure}[h] 
  \begin{center}
    \includegraphics[width=.7\textwidth]{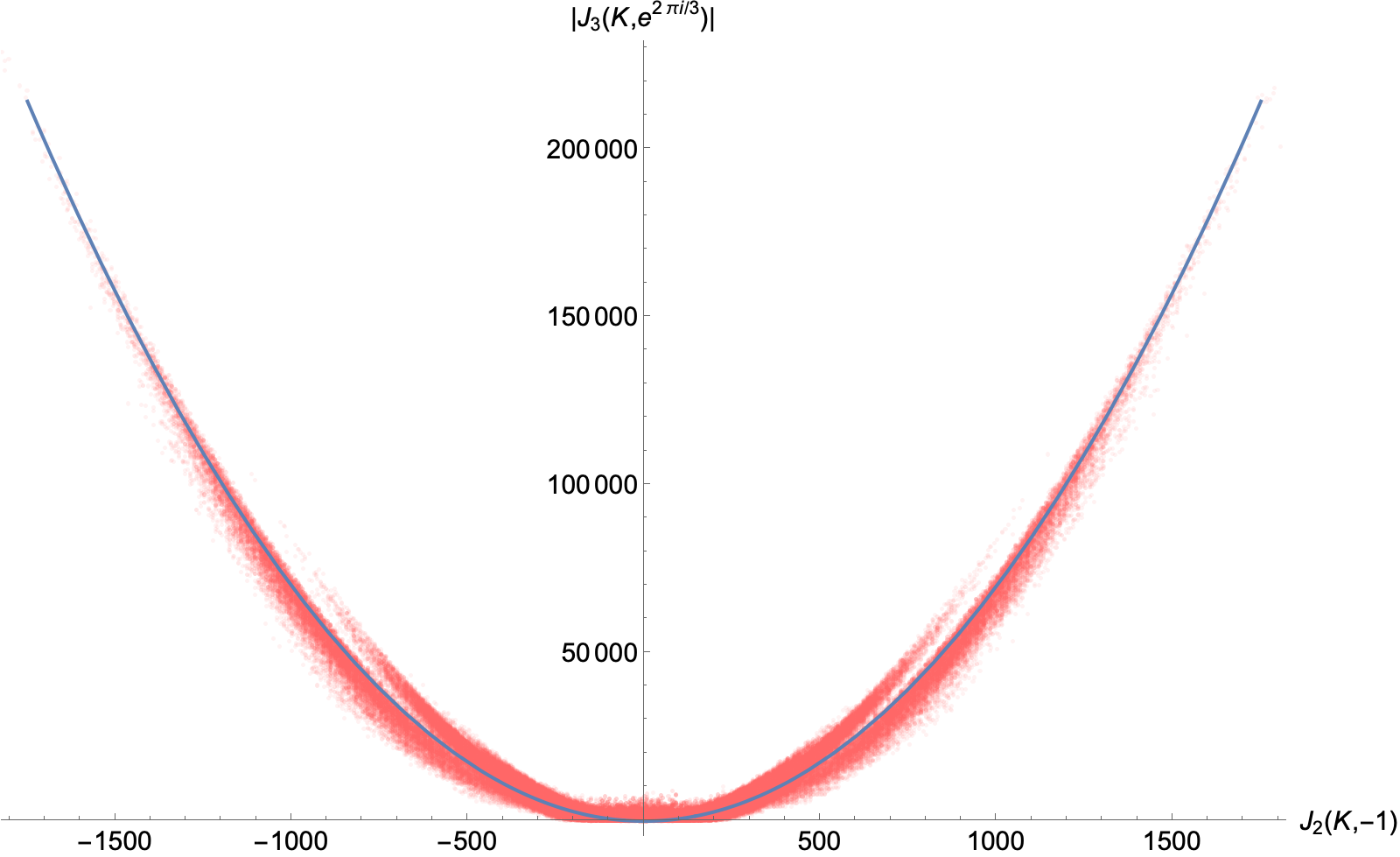}
  \end{center}
  \caption{\textsf{Plot of evaluations of $J_3$ at $e^{2\pi i/3}$ vs.\ evaluations of $J_2$ at $e^{-i\pi}=-1$. In blue is plotted the curve $x^2/10^{1.155}$, which fits the data nicely. This is in accord with a similar plot with a smaller dataset in~\cite{Craven:2020bdz}.}}  \label{fig:quadratic-relation}
\end{figure}

As was observed in~\cite{Craven:2020bdz}, we also find a quadratic relation between the evaluation of $J_2$ at $e^{i\pi}=-1$ and of $J_3$ at $e^{2\pi i/3}$.
These are the phases implicated by the statement of the volume conjecture,~\eref{eq:vc}.
This is plotted in Figure~\ref{fig:quadratic-relation}.
We find that the curve $10^{1.155}y=x^2$ fits the data nicely, as also shown in the figure.

\section{Results from training neural networks}\label{sec:ml}

We now present results from a variety of supervised machine learning (ML) experiments on the fundamental and adjoint Jones polynomials, with the aim of determining correlations between the polynomials and the hyperbolic volumes of knots.
We then use the results to arrive at a guess for the best phase for use in the volume conjecture and check the predictions from this phase.

\subsection{Using polynomials}
A review of neural networks and their utility in mathematics is~\cite{williamson2024deep}.
Here, we present results obtained by using fully connected feedforward neural networks (NNs) with five hidden layers.
The architecture of this deep neural network can be represented as%\\ \mbox{$(\text{Input},100,\text{ Ramp}, 300,\text{ Ramp},300,\text{ Ramp},150,\text{ Ramp},75,\text{ Ramp},\text{Output})$}, 
\[
    (\text{Input},100,\text{ReLU}, 300,\text{ReLU},300,\text{ReLU},150,\text{ReLU},75,\text{ReLU},\text{Output}) \, ,
\]
where the numbers are the numbers of neurons in each of the hidden layers and ReLU is the $\text{max}(0,x)$ non-linearity.
We use the ADAM optimizer with an adjustable learning rate.
The results are insensitive to minor tweaks to the neural network architecture.
The input to the networks are vectors encoding the content of the Jones polynomials.
We formed the vectors by taking the first two components of a vector to be the minimum and maximum degrees of the corresponding polynomial and the next components of the vector to be coefficients of all monomials from lowest to highest degree.
We then padded the vectors with zeros on the right to make them of uniform length.
For $J_2$, we pad the vectors to make them of length $19$.
For $J_3$, we pad these vectors to have length $48$.
We also present results obtained by training deep neural networks on vectors formed from joining the vectors for $J_2$ and $J_3$ of a given knot.
In this case we pad $J_2$ vectors to have the same length as $J_3$ vectors, resulting in vectors of length $96$ in this case.

\begin{figure}[h] 
  \begin{center}
    \includegraphics[width=.7\textwidth]{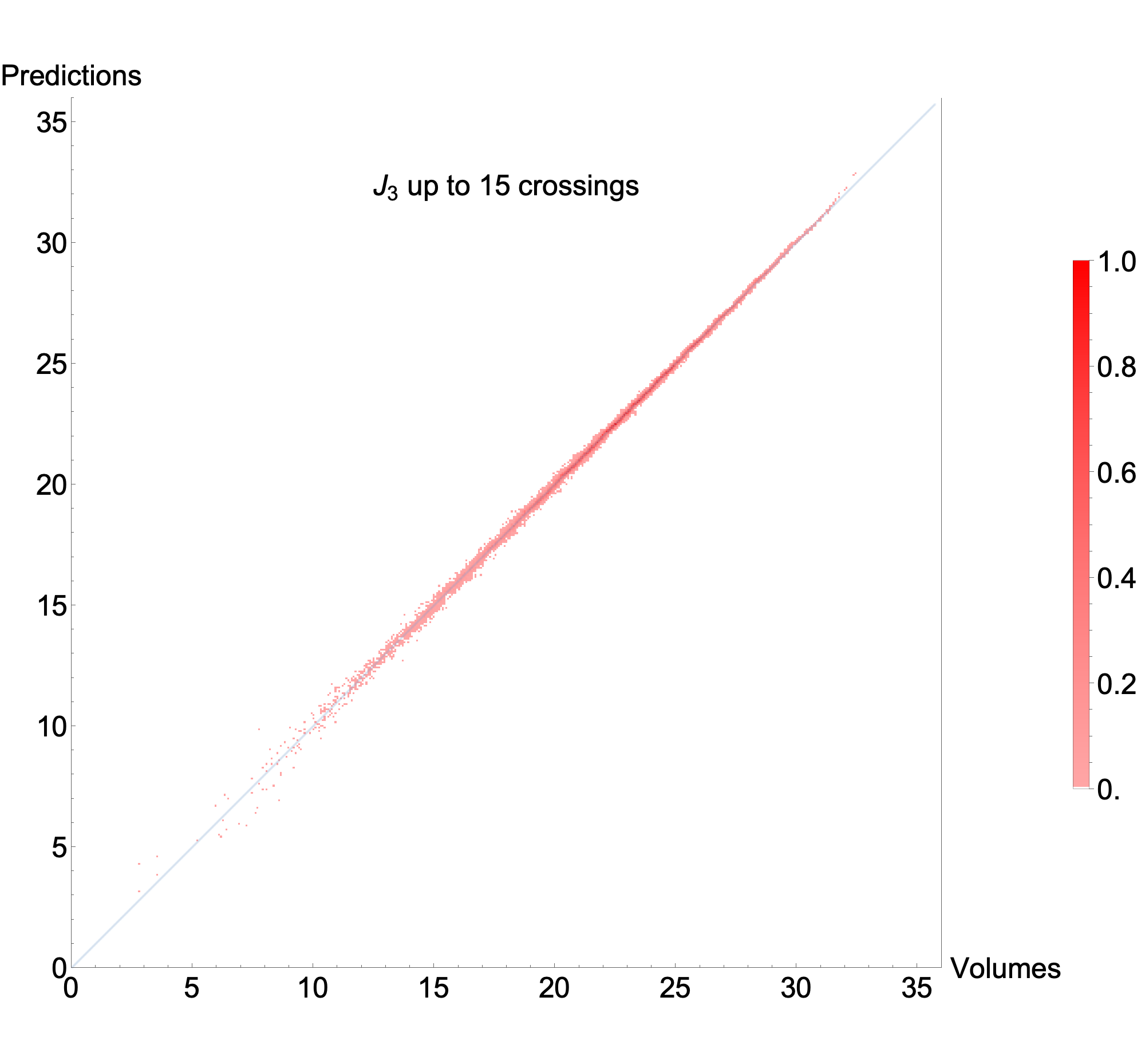}
  \end{center}
  \caption{{\textsf{Predicted vs.\ actual volumes using a neural network trained on $J_3$ for knots with up to $15$ crossings. Also included are the actual volumes in blue in the background. The color signifies the density of points, as depicted by the legend bar. The mean relative error of predictions (on the test set) by the neural network used to generate this plot was $0.66\%$. Note that here and in all similar plots below, the intensity of points signifies the density of knots in each pixel of the plot.}}}  \label{fig:J3poly131415preds}
\end{figure}

Using $J_3$ polynomials of knots up to $15$ crossings, the mean relative error\footnote{All the reported relative errors are for performance over test samples.} over three independent runs of training (with $100$ epochs each) was $0.40\%$.
In the experiments, we used $75\%$ of the knots for training, $10\%$ for validation, and the remaining $15\%$ for testing the trained network.
When the $J_3$ polynomials of knots with only $14$ and $15$ crossings were used, the mean relative error over three runs was $0.37\%$.
Using only knots of $15$ crossings, the mean relative error over three runs was $0.35\%$.
A plot of predictions from a neural network trained on $J_3$ polynomials is included in Figure~\ref{fig:J3poly131415preds}.

\begin{figure}[h] 
  \begin{center}
    \includegraphics[width=.7\textwidth]{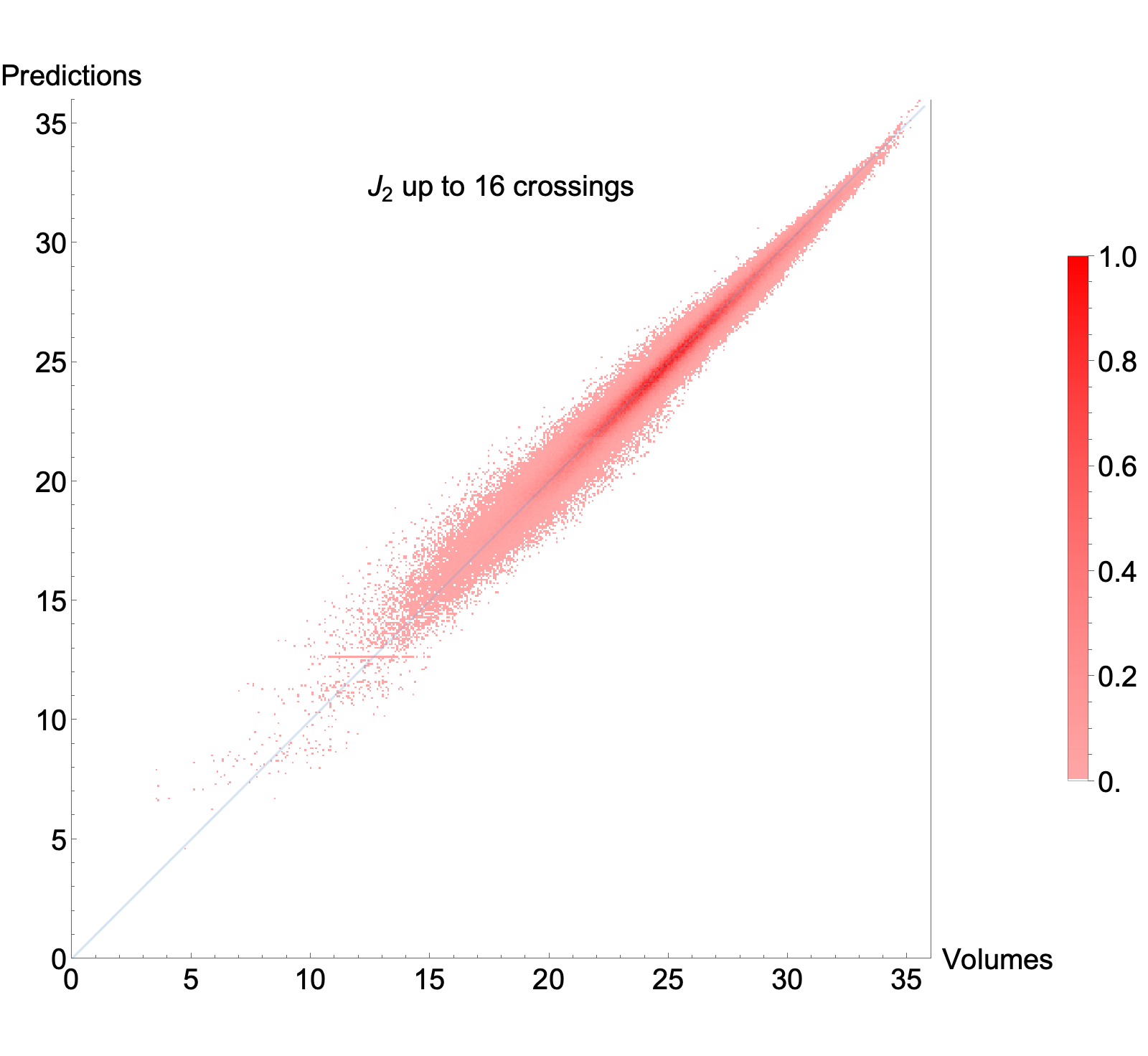}
  \end{center}
  \caption{{\textsf{Predicted vs.\ actual volumes using a neural network trained on $J_2$ for knots with up to $16$ crossings. Also included are the actual volumes in blue. The mean relative error of predictions (on the test set) by the neural network used to generate this plot was $1.64\%$. This is consistent with the results reported in~\cite{Jejjala:2019kio, Craven:2020bdz}, where only $10\%$ of the dataset was used for training.}}}  \label{fig:J2poly131415preds}
\end{figure}

For machine learning using $J_2$ polynomials of the $1{,} 701{,} 903$ hyperbolic knots up to $16$ crossings, the mean relative error over three independent runs of training was $1.65\%$.
For comparison, a plot of predictions from a neural network trained on $J_2$ polynomials of all hyperbolic knots up to $16$ crossings is included in Figure~\ref{fig:J2poly131415preds} (\textit{cf.},~\cite{Jejjala:2019kio, Craven:2020bdz}).
With training on the hyperbolic knots up to $15$ crossings, the mean error over three training runs was $1.86\%$.
When restricted to the $177{,}316$ knots for which we know the $J_3$ polynomials, this error does not change significantly.
When the $J_2$ polynomials of knots with only $14$ and $15$ crossings were used, the mean relative error over three runs was $1.55\%$.
When knots of only $15$ crossings were used, the mean relative error over three runs was $1.48\%$.
When knots of only $16$ crossings were used, the mean relative error over three runs was $1.44\%$.

The results reported above seem to suggest that the correlations between Jones polynomials and the volumes of knots improves as the number of crossings of knots increases.
To check if this is indeed the case, we restricted to random samples of $100,000$ knots of $15$ and $16$ crossings each and trained neural networks on the $J_2$ of each set separately.
Over three training runs, the mean error from the neural networks trained on $15$ and $16$ crossing knots was $1.695\%$ and $1.697\%$ respectively.
It would thus seem that it is the greater number of knots at higher crossings that results in better performance for the neural networks.

\begin{figure}[h] 
  \begin{center}
    \includegraphics[width=.7\textwidth]{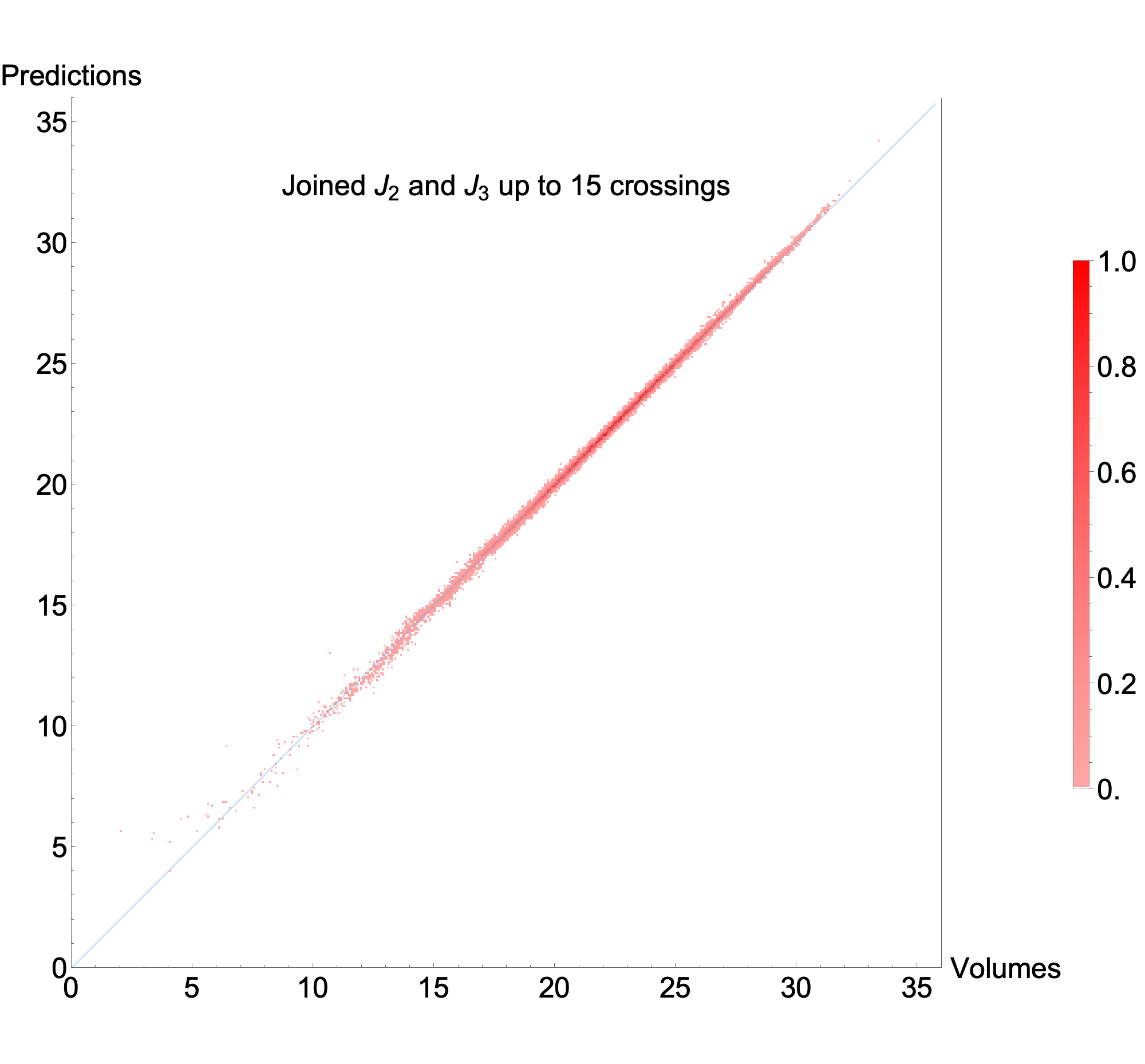}
  \end{center}
  \caption{{\textsf{Predicted vs.\ actual volumes using a neural network trained on the joined vectors of coefficients formed from $J_2$ and $J_3$ polynomials, for knots with up to $15$ crossings. Also included are the actual volumes in blue. The mean relative error of predictions (on the test set) by the neural network used to generate this plot was $0.54\%$.}}}  \label{fig:J2J3poly131415preds}
\end{figure}

For machine learning with $J_2$ and $J_3$ combined into one vector of coefficients, the mean error over three runs of training was $0.51\%$.
This is shown in the right plot in  Figure~\ref{fig:J2J3poly131415preds}.

To better compare neural networks trained on different colored polynomials, we trained and tested neural networks on the same dataset.
The mean relative errors of the trained networks using $J_2$, $J_3$, and joined $J_2$ and $J_3$ polynomials was respectively $1.85\%$, $0.62\%$, and $0.40\%$.

To summarize, we find that the $J_3$ polynomials predict the volumes with significantly more accuracy than $J_2$ polynomials.
The performance of joined $J_2$ and $J_3$ polynomials is roughly the same as that of just $J_3$ polynomials. 
Thus, despite the non-monotonicity reported in~\cite{Garoufalidis:2004} based on the study of a specific knot, it would appear that statistically the approach to the volume conjecture using a large dataset might be monotonic.

\subsection{Using evaluations}

\begin{figure}[h]
    \centering
    \includegraphics[width=.6\textwidth]{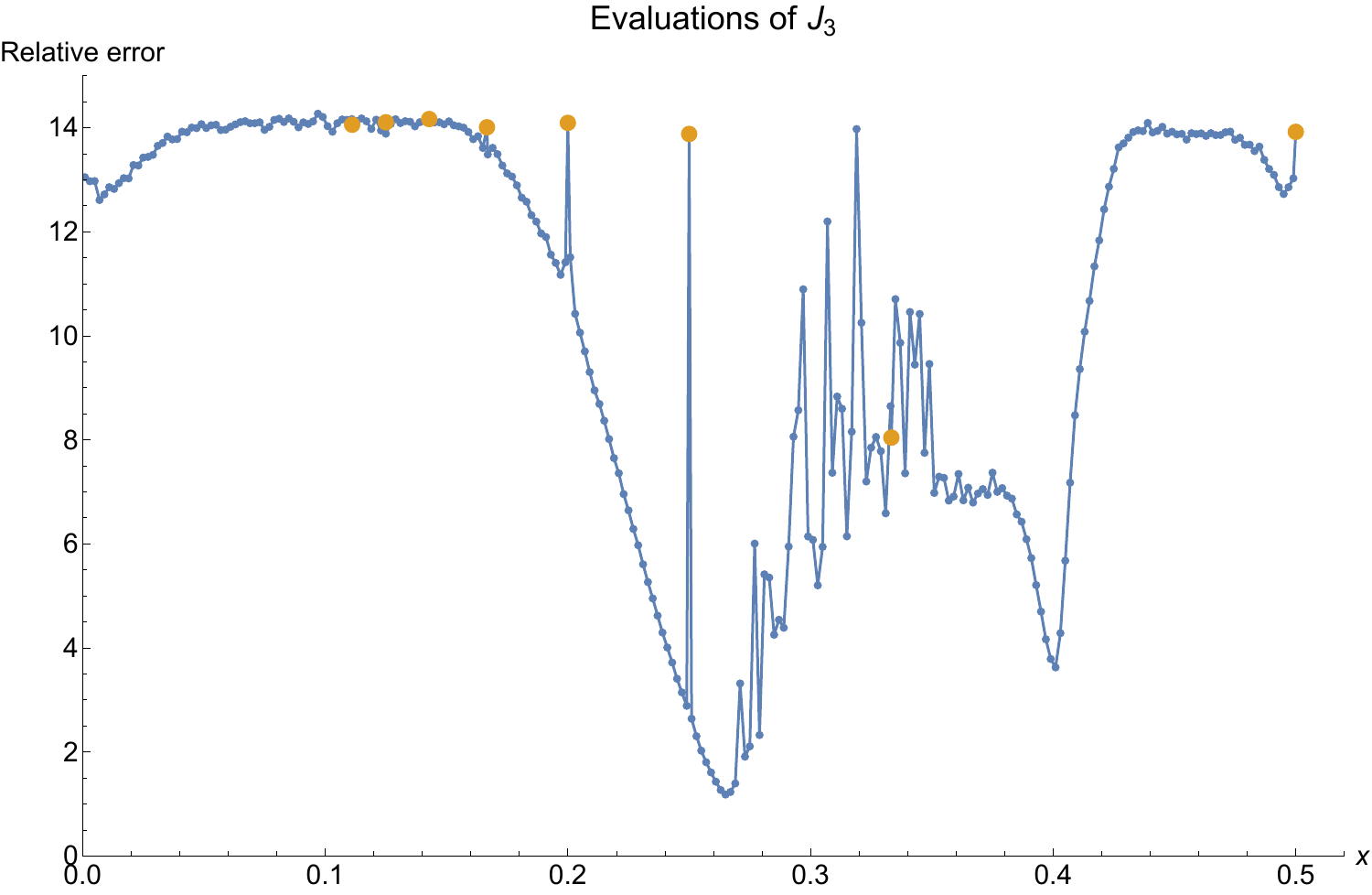} 
    \caption{\textsf{Relative error of predictions of neural network with input the evaluations $J_3(K;e^{2\pi i x})$. This is over three runs of training with $25$ epochs. The orange points have $x=(k+2)^{-1}$ for integer $k\in[0,7]$.}}
    \label{fig:J3-phases-rel-err}
\end{figure}

In Figure~\ref{fig:J3-phases-rel-err}, we plot the relative error of the predictions of the neural network with respect to the actual volume as a function of the phase at which the $J_3$ polynomials are evaluated.
The network was given as inputs the real and imaginary parts of the evaluation of the polynomial at the phases.
We checked that giving the modulus and argument of the evaluation as input does not change the results significantly.
The best phase is seen to be $x\sim0.27$.
We also performed similar experiments with input the sum of $|J_2(K;e^{2\pi ix})|+|J_3(K;e^{2\pi ix})|$, and found a qualitatively similar plot.
This is expected since the evaluations of $J_3(K;q)$ tend to have larger magnitude than the evaluations of $J_2(K;q)$.
The existence of a good minimum is again explainable via the analytically continued Chern--Simons theory~\cite{witten2011analytic} and tracks the discussion in~\cite{Craven:2020bdz}.
It is not \textit{a priori} apparent why there should be two minima in the phase vs.\ error plot for the $3$-color Jones polynomial compared to only one minimum in the $2$-colored case --- \textit{cf.}, Figure~7 in~\cite{Craven:2020bdz}.

\begin{figure}[h] 
  \begin{center}
    \includegraphics[width=.6\textwidth]{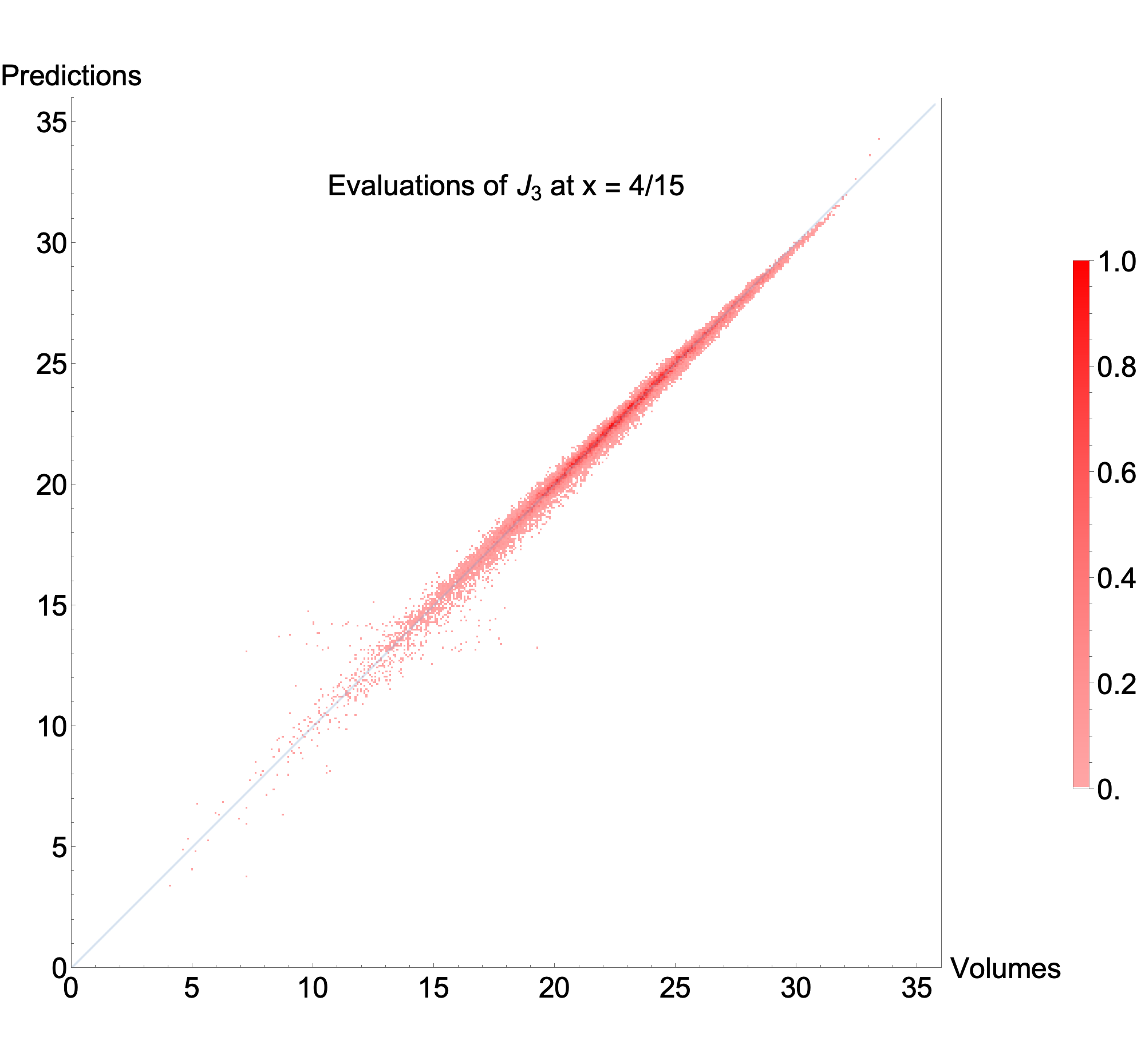}
  \end{center}
  \caption{\textsf{Predicted vs.\ actual volumes using a neural network trained on $J_3(K;e^{2\pi ix})$ with $x=4/15$. Also included are the actual volumes in blue. The mean relative error of predictions by the neural network used to generate this plot was $1.14\%$.}}  \label{fig:J3-evals-4by15}
\end{figure}

In Figure~\ref{fig:J3-evals-4by15}, we show the predictions for volume obtained by using a neural network trained on the evaluations $J_3(K;e^{2\pi ix})$ at $x=4/15\approx 0.26667$, which is the best phase in Figure~\ref{fig:J3-phases-rel-err}.
The mean relative error of the predictions of this neural network is $0.90\%$.
Thus evaluations of $J_3$ at this phase perform only slightly worse than the data containing the full polynomials.
Using the absolute value of evaluation of polynomials at $x=4/15$ to train neural networks, we find that the mean error over three runs increases slightly to $1.33\%$.

If we use evaluations of $J_3$ at $e^{2\pi i\,2/5}$, which is the other minimum in Figure~\ref{fig:J3-phases-rel-err}, to train a neural network for predicting the volume, the mean relative error on the test set is $3.60\%$.
When using the sum of evaluations at $e^{8\pi i/15}$ and $e^{4\pi i/5}$, the mean relative error improves, and becomes $1.61\%$.

\begin{figure}[h] 
	\begin{center}
    \includegraphics[width=.5\textwidth]{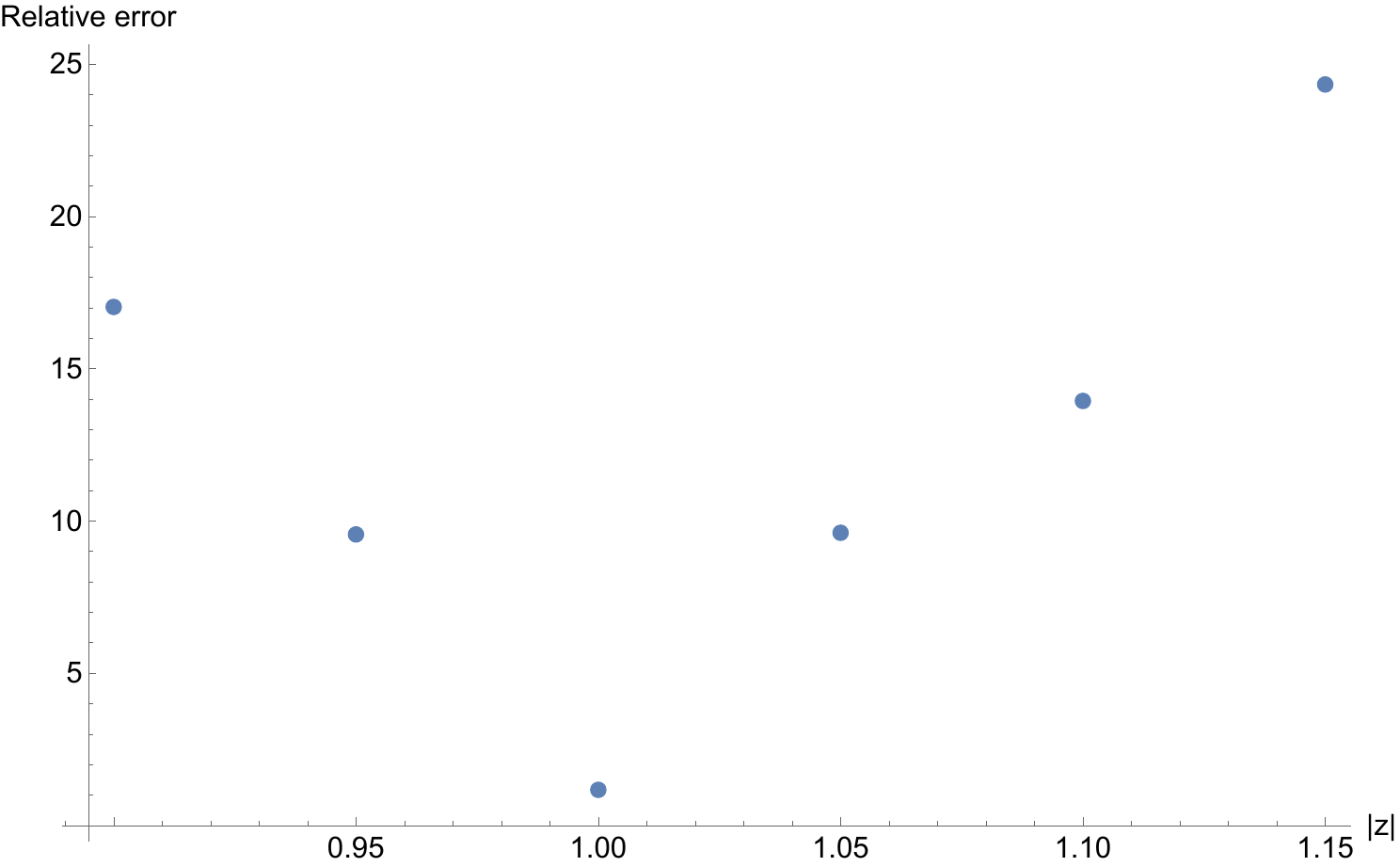}
    \end{center}
  \caption{\textsf{Plots of relative errors of neural networks trained on evaluations at points $z=|z|e^{2\pi ix}$ with different absolute values with fixed phase at $x=4/15$.}}  \label{fig:off-circle}
\end{figure}

Finally, in Figure~\ref{fig:off-circle}, we show the errors of predictions for neural networks trained using evaluations off the unit circle.
The error increases monotonically as the radial distance from the unit circle is increased on either side for the three phases that we have checked.


\subsection{Symbolic formula for the volume}

\begin{figure}[t]
    \begin{center}
        \includegraphics[width=.5\textwidth]{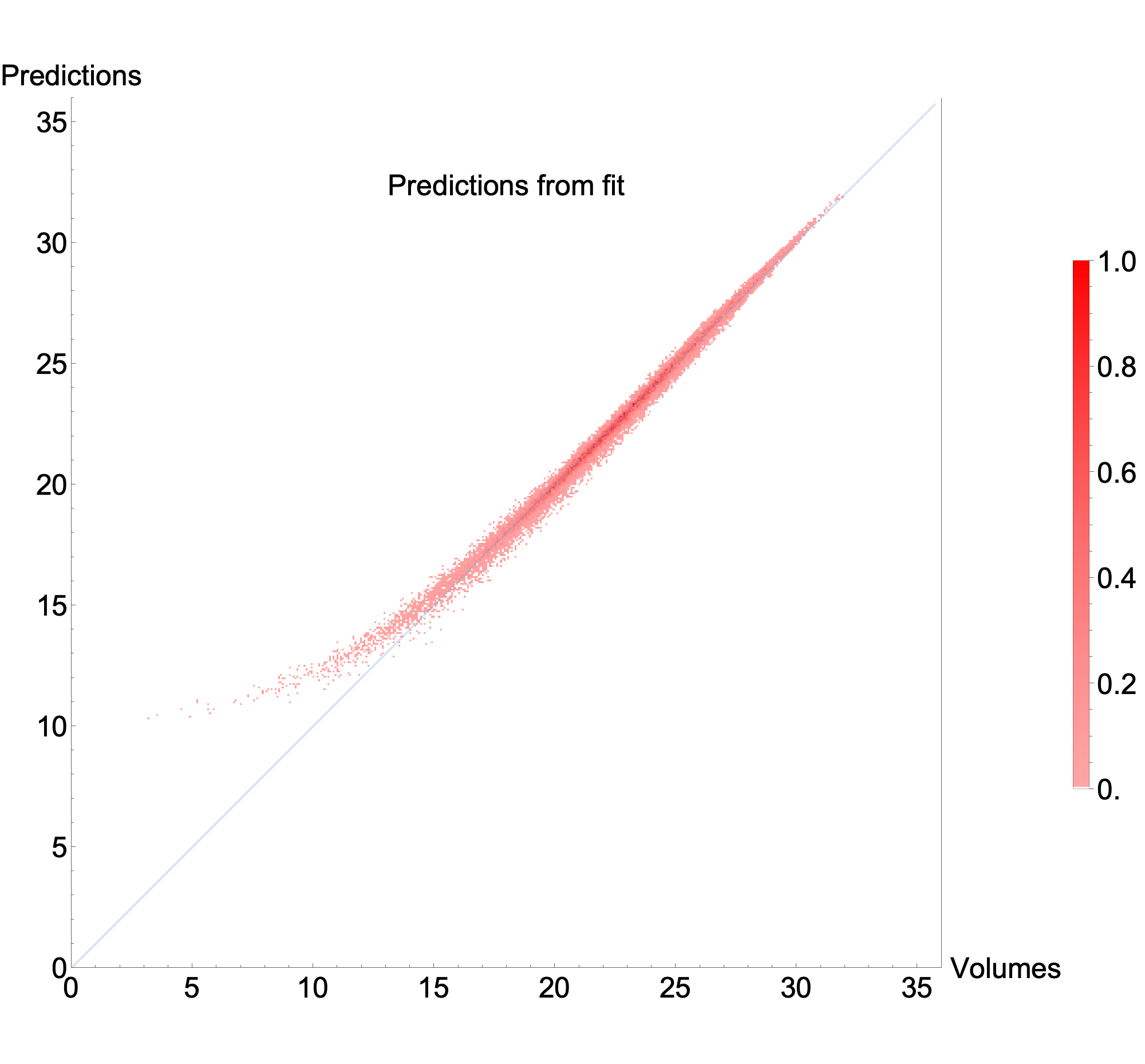}
    \end{center}
    \caption{\textsf{Predictions from the symbolic expression~(\ref{eq-regression-result}).}}  \label{fig:nlm-fit}
\end{figure}

Using the Mathematica function \texttt{NonlinearModelFit} with a guess
\begin{equation}
    a \log(b|J_3(K;e^{8\pi i/15})|+c)+d
\end{equation}
for the functional form and evaluations at $x=4/15$ as the input data to the function, we get the following function for volume,\footnote{
Using the Julia/Python package \texttt{PySR} for symbolic regression with $\exp$ and $\log$ as guesses for the functions involved, we get a similar formula as the one obtained using Mathematica.}
\begin{equation}			\label{eq-regression-result}
	\text{vol}_K = 3.25\log(|J_3(K;e^{8\pi i/15})|+36.97) - 1.72 \,.
\end{equation}
The mean relative error for predictions of volumes using this formula is $1.21\%$. As further measure of the quality of the fit, we use the coefficient of determination, usually denoted $R^2$, with $R^2=1$ signifying an exact formula.
This formula gives an $R^2$ value of $0.999625$ and is thus a very good fit.
Predictions from this formula are plotted in Figure~\ref{fig:nlm-fit} for a random sample of $25{,} 000$ knots.

\section{Toward an improved volume conjecture}\label{sec:newvc}

In~\cite{Craven:2020bdz}, the phase at which the evaluation of the $2$-colored Jones polynomials was found to give a best fit to the volume was $e^{3\pi i/4}$.
We have presented evidence that $e^{8\pi i/15}$ is the corresponding best phase to evaluate $3$-colored Jones polynomials to train neural networks to predict the volumes of knot complements of hyperbolic knots.
These correspond, respectively, to fractional Chern--Simons levels $k = 2/3$ and $7/4$ and $\gamma = 3/2$ and $8/7$.
Using these two data points, we make the guess that for $n$-colored Jones polynomials, the best phase to use is
\begin{equation}        \label{eq:best-phase-guess-2}
    q(n) = \exp\left( 2\pi i \frac{n+1}{n(n+2)} \right).
\end{equation}
Correspondingly, we find that the fractional Chern--Simons level and $\gamma=(n-1)/k$ are
\begin{equation}
    k(n) = \frac{n^2-2}{n+1} \,, \qquad
    \gamma(n) = \frac{n^2-1}{n^2-2} \,.
\end{equation}
These formul\ae\ have the expected asymptotics.
In the large-$n$ limit,
\begin{equation}
q(n) \sim \exp \left( 2\pi i \left( \frac{1}{n}-\frac{1}{n^2} \right) \right) \,, \qquad
k(n) \sim (n-1)-\frac{1}{n} \,, \qquad
\gamma(n) \sim 1+\frac{1}{n^2} \,.
\end{equation}
This tells us the relations are in accord with the volume conjecture,~\eref{eq:vc}, which posits that at large-$n$, the phase at which the $n$-colored Jones polynomial is evaluated is the primitive $n$-th root of unity.
The machine learning experiments suggest that rather than evaluating at the phase $\omega_n = 2\pi i/n$, we should dress this as
\begin{equation}
2\pi i\; x(n) = \frac{n+1}{n+2}\; \omega_n \,.
\end{equation}
Thus, we propose that~\eref{eq:vc} be modified to read
\begin{equation}\label{eq:vc-new}
	\lim_{n\to\I} \frac{\log|J_n(K; q(n))|}n = \frac1{2\pi} V(S^3\backslash K) \,.
\end{equation}

We test this hypothesis using exact expressions for $n$-colored Jones polynomials available for the two hyperbolic knots $4_1$ and $K_0$.

The $n$-colored Jones polynomials for $4_1$ are~\cite{Le:2000, Habiro:2000, Murakami:2022}
\begin{equation}
    J_n(4_1;q) = \sum_{k=0}^{n-1} q^{-nk} \prod_{l=1}^k (1-q^{n+l}) (1-q^{n-l}) \,.
\end{equation}
From Snappy, the volume of $4_1$ is $2.02988321282$.

The $(n+1)$-colored Jones polynomials of $K_0$ are~\cite{Garoufalidis:2004}
\begin{align}
  J_{n+1}(K_0;q) =&\ \frac1{[n+1]}\sum_{k=0,2}^{2n}\sum_{l=|n-k|,2}^{n+k}\sum_z(-1)^{\frac k2+z}q^{-\frac38(2k+k^2)+\frac78(2l+l^2)-\frac{51}8(2n+n^2)}\frac{[k+1][l+1]}{[n+\frac k2+1]!}   \nn \\
  & \quad \times \begin{bmatrix} \frac{k+l-n}2 \\ \frac{n+2k+l}2 - z  \end{bmatrix} \begin{bmatrix} \frac{n+l-k}2 \\ \frac{3n+l}2 - z  \end{bmatrix} \begin{bmatrix} \frac{n+k-l}2 \\ n+k - z  \end{bmatrix} \lB\frac k2 \rB!^2 \frac{\lB n-\frac k2 \rB!}{\lB z-\frac{n+k+l}2 \rB!}\frac{[z+1]!}{\lB \frac{n+k+l}2+1 \rB!} \ ,   
\end{align}
where $\sum_{k=a,2}^b$ means summation with step $2$ from $a$ to $b$, and the summation over $z$ has the limits
\begin{equation}
    \max\lb n+\frac k2,\frac{n+k+l}2 \rb \leq z\leq \min\lb \frac{n+2k+l}2,\frac{3n+l}2, n+l \rb \,.
\end{equation}
We note that the volume of $K_0$ is $3.474247$.

\begin{figure}
	\centering
    \subfloat[$4_1$]{%
        \includegraphics[width=.45\textwidth]{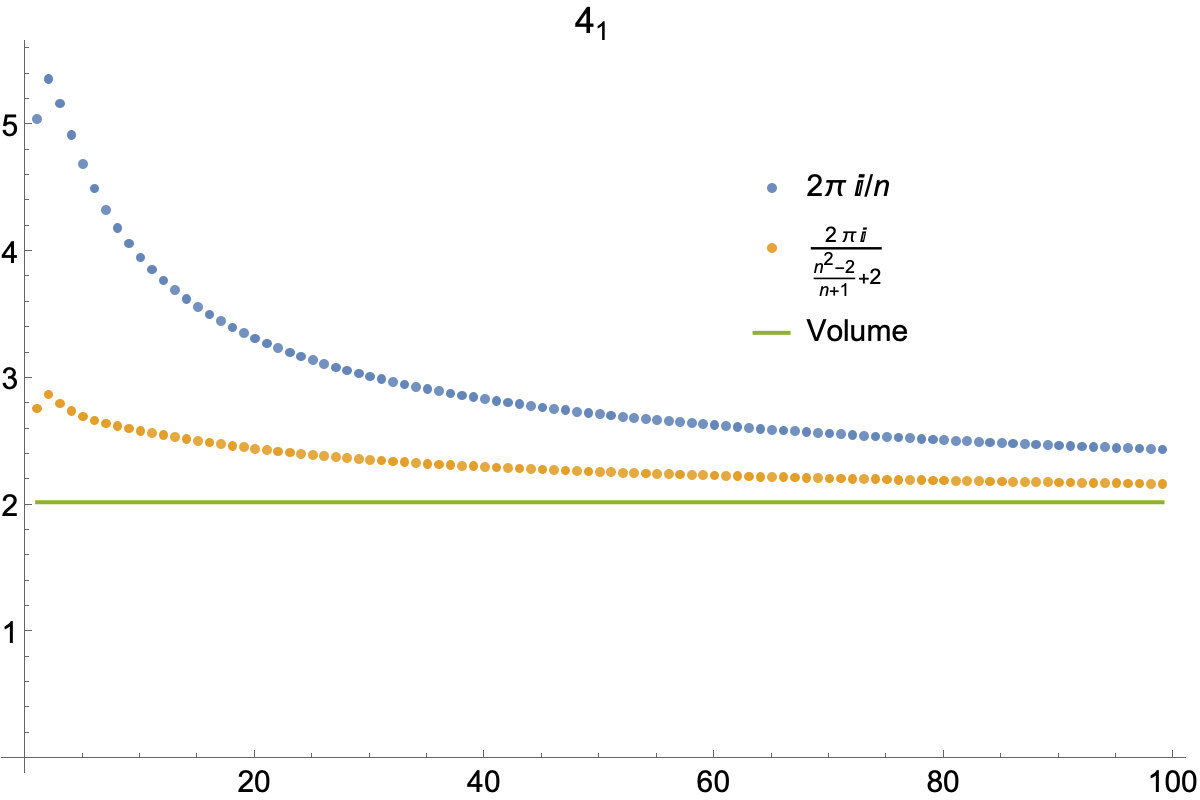}	\label{fig:41}%
    }
    	\subfloat[$K_0$]{%
    	\includegraphics[width=.45\textwidth]{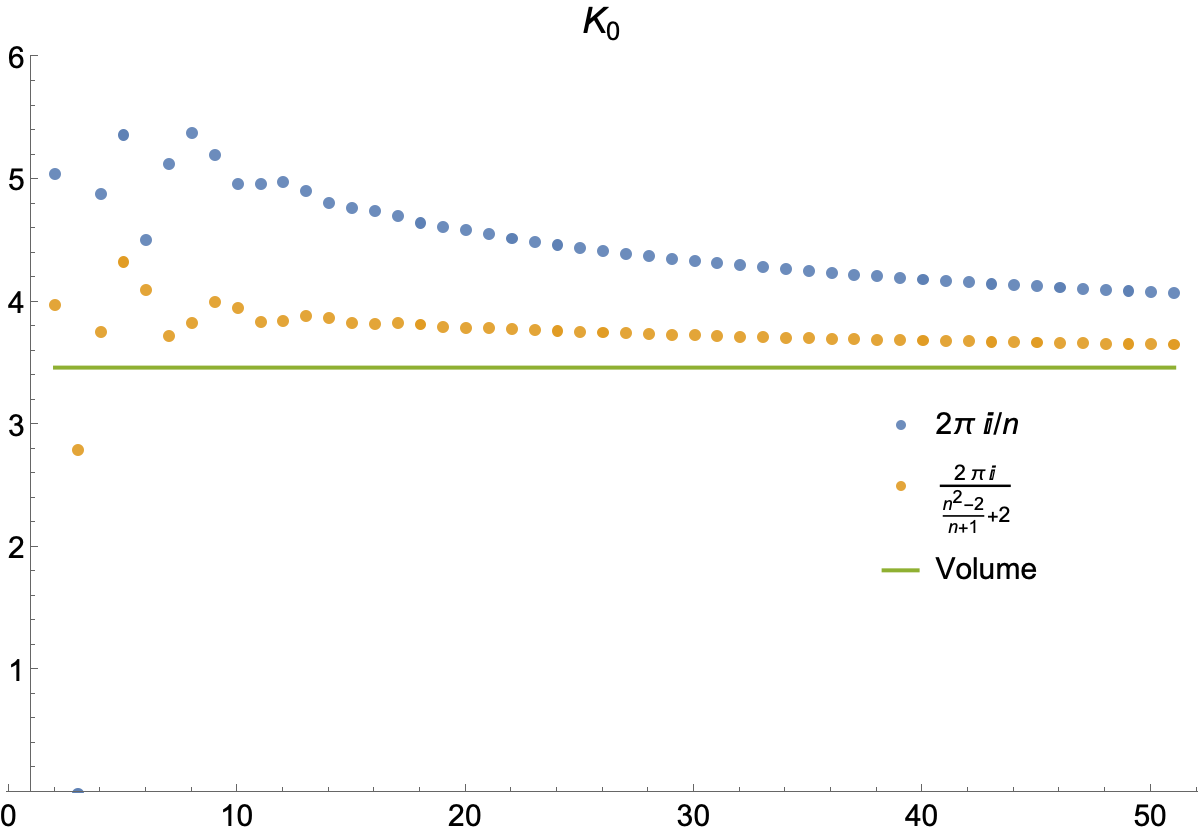}     \label{fig:K0}
	} 
    \caption{\textsf{Comparison of predictions using the phase $\omega_n$ from~(\ref{eq:vc}) and using the phase from~(\ref{eq:best-phase-guess}), along with the actual volumes of the knots $4_1$ (left) and $K_0$ (right).}}
    \label{fig:K0-41}
\end{figure}

We use the knowledge of these exact polynomials to test whether the conjectured phase $q(n)$ in~(\ref{eq:best-phase-guess-2}) results in improvements upon the prediction of the volume at finite $n$.
In Figure~\ref{fig:K0-41}, along with the volumes of the knots and the left hand side of~\eref{eq:vc} (the volume conjecture), we plot the values of
\begin{equation}\label{eq:vc2}
	v(n) = \frac{2\pi\log|J_n(K;q(n))|}n \,.
\end{equation}
We see that the phase $q(n)$ performs much better for both the figure-eight knot and the knot $K_0$.
We can additionally envision a subleading shift in the denominator to achieve a better convergence, but we have not systematically investigated this possibility.
The convergence of $v(n)$ to the volume is still not monotonic.
Thus, it is impressive that the neural network does so much better using the adjoint Jones polynomial, and this should be understood.
It would also be worthwhile to test whether our conjectured phase continues to perform well for other knots at higher colors.
This requires a database of higher colored Jones polynomials.
Work on this is underway, and we hope to report on progress soon.

\section*{Acknowledgements}
The content of this paper was presented at a workshop on the ``Generalized Volume Conjecture'' at IIT Bombay in November 2024 and the String Data meeting at the Yukawa Institute for Theoretical Physics in Kyoto in December 2024.
We thank the organizers and participants at these meetings for their generous feedback.
We are especially grateful to Sergei Gukov and Hisham Sati for discussions.
We thank Aditya Dwivedi for conversations about improving the databases of colored Jones polynomials.
MH is supported by the National Science Foundation (DMS-2213295), and thanks the Dublin Institute for Advanced Studies and the Max Planck Institute for Mathematics for hosting him during work on this paper.
VJ is supported by the South African Research Chairs Initiative
of the Department of Science and Innovation and the National Research Foundation.
P.~Ramadevi would like to thank the IITB-IOE funding:  `Seed funding for Collaboration and Partnership Projects (SCPP) scheme - Phase III'.
P.~Roy is supported by an NRF Freestanding Postdoctoral Fellowship. P.~Roy would like to thank Chennai Mathematical Institute for hospitality while this work was being completed.
VKS is supported by the ‘Tamkeen under the NYU Abu Dhabi Research Institute grant CG008 and ASPIRE Abu Dhabi under Project AARE20-336. This work was supported in part through the High-Performance Computing resources at New York University Abu Dhabi.

\appendix
\section{Knot polynomials from vertex models}\label{sec:app}

There are various methods to compute polynomial invariants for knots and links~\cite{Witten:1988hf,KirResh,RT1,RT2,Itoyama:2012qt,Kaul:1992rs,Devi:1993ue,RGK,DMMMRSS:2017,DMMMRSS:2018}. For our calculations of the adjoint Jones polynomials of knots with $14$ and $15$ crossings, we adopted a vertex model approach, which we found to be significantly faster than the \texttt{KnotTheory} package for Mathematica available at KnotAtlas~\cite{KnotAtlas}. We now briefly review our approach to constructing Jones polynomials using braid group representations obtained via vertex models. This material is well known in the literature, and we follow the review~\cite{Wadati:1989ud}. 

Vertex models on square lattices are statistical mechanical systems where each vertex has four associated edges, and the states on the four edges determine the Boltzmann weights of the vertices. Each vertex can be thought of as describing two-to-two scattering with particles of various charges. One can study vertex models where the edge variables are allowed to take $N$ possible states. Thinking of these states as ``spin'' charges, and demanding appropriate conservation, one obtains what are called $N$-state vertex models with spin $s=(N-1)/2$. These models can all be shown to be integrable. 

The braid group over $n$ strands, denoted $B_n$, has $n-1$ generators $b_i$, $i=1,\ldots,n-1$, subject to the relations
\begin{eqnarray}
    b_i\, b_{i+1} \, b_i = b_{i+1} \,b_i \, b_{i+1} & \text{ for } & 1\leq i \leq n-2 \,, \text{ and } \nonumber \\
    b_i\, b_j = b_j\, b_i & \text{ for } & |i-j|\geq 2 \,.
\end{eqnarray}
Any $n$-braid is expressed as a word in $B_n$, \textit{e.g.}, $b_1^{-1}b_2b_3b_2^{-1}$, up to the above relations. Given a braid, one forms an oriented link by taking its closure, \textit{i.e.}, by identifying opposite ends of the braid. Conversely, any oriented link is represented (non-uniquely) by some closed braid. Markov's theorem states that two closed braids represent the same ambient isotopy class of links if and only if the braids can be transformed into one another by a sequence of ``Markov moves'' of type I and II, also called conjugation and stabilization. The moves are
\begin{eqnarray}
	\text{(I)} &\quad& AB\to BA \,, \text{ for } A,B\in B_n \,, \\
	\text{(II)} &\quad& A\to A\,b_n^{\pm1}\,, \text{ for } A\in B_n, \ b_n\in B_{n+1} \,.
\end{eqnarray}
One can then obtain invariant polynomials for links by constructing a representation of the braid group, and then defining Markov move invariant polynomials using the representation.

\subsection{Akutsu--Wadati formula}

Using the Boltzmann weights of vertex models, one can construct the so-called Yang--Baxter operators satisfying Yang--Baxter equations. Upon ``asymmetrizing'' and taking appropriate limits, these Yang--Baxter operators furnish representations $G_i$ of braid operators $b_i$, of the form 
\begin{equation}
	G_i = \mathbb{I}_1\otimes\mathbb{I}_2\otimes\cdots\otimes\mathbb{I}_{i-1}\otimes R_{i,i+1}\otimes\mathbb{I}_{i+1}\otimes\cdots\otimes\mathbb{I}_{n}\,,
\end{equation}
where the matrix $R$ acts on two strands of the braid, and $\mathbb{I}_j$ denotes the appropriate identity matrix acting on the strand $j$. We note that the matrix $R$, called the (braided) $R$-matrix, is fully determined by the specification of a given integrable statistical model~\cite{Baxter:1982zz}. 

Let $G_i$ be such a representation of $b_i\in B_n$, and let $\a(\cdot)$ denote a link polynomial. From the above discussion, we know that $\a(\cdot)$ must satisfy the following conditions:
\begin{eqnarray}
	\text{(I)} &\quad& \a(AB) = \a(BA) \,, \text{ for } A,B\in B_n \,, \\
	\text{(II)} &\quad& \a(A)= \a(AG_n)=\a(AG_n^{-1}) \,, \text{ for } A\in B_n,G_n\in B_{n+1} \,.
\end{eqnarray}
Such an $\a(\cdot)$ can be constructed if we can construct a Markov trace, $\f(\cdot)$, satisfying
\begin{eqnarray}
	\text{(I)} &\quad& \f(AB) = \f(BA) \,, \text{ for } A,B\in B_n \,, \\
	\text{(II)} &\quad& \f(A G_n) = \t\f(A) \f(AG_n^{-1})=\bar\t\f(A) \,, \text{ for } A\in B_n\,,\,G_n\in B_{n+1} \,,
\end{eqnarray}
where $\t=\f(G_i),\ \bar\t=\f(G_i^{-1})$ for any $i$. In terms of the Markov trace $\f(\cdot)$, the link polynomial is given by the following (Akutsu--Wadati) formula,
\begin{equation}
	\a(A) = (\t\bar\t)^{-(n-1)/2}\lb\frac{\bar\t}\t\rb^{e(A)/2}\f(A) \,, \quad A\in B_n \,,
\end{equation}
where $e(A)$ is the sum of the exponents of $b_i$ in $A$. For example, for the braid word $A=b_1 b_2^{-1} b_1 b_2^{-1}$, $e(A)=0$.
For $N$-state vertex models, we have~\cite{Akutsu:1987dz,Akutsu:1987qs}
\begin{equation}
	\t = \frac1{1+q+\cdots+q^{N-1}}\,, \qquad \bar\t = \frac{q^{N-1}}{1+q+\cdots+q^{N-1}}\, .
\end{equation}
The Markov trace is given explicitly by
\begin{equation}
	\phi(A) = \text{Tr}(H \cdot A) \,,
\end{equation}
where $H$ is the tensor product of $n$ matrices $h$ of size $N\times N$,
\begin{equation}
	H = \underbrace{h \otimes h \otimes \dots \otimes h}_n \,.
\end{equation}
The matrix $h$ is diagonal, and is given for $N$-state vertex models by
\begin{equation}
	h = \t\,\text{diag}\,(1,q,\cdots,q^{N-1}) \,.
\end{equation}
The only remaining element is then the braided $R$-matrix, which needs to be computed just once for each $N$-state vertex model. 

It can be shown that the link invariant polynomials that we have constructed using the $N$-state vertex model and the Akutsu--Wadati formula correspond, for each $N$, to the $N$-dimensional Jones polynomials of links. Below, we give a few example calculations of Jones polynomials for the trefoil and figure-eight knots.

Before proceeding we make a remark about the complexity of this method. The matrix $H$ obtained from an $N$-state vertex model is of size $N^n\times N^n$ for a braid word of $n$ strands. Since the $R$-matrix acts on the Hilbert space of two ``particles'', it is a matrix of size $N^2\times N^2$. The braid generators $b_i$ will still have size $N^n\times N^n$. Because of this, it quickly becomes prohibitively expensive to use this method to calculate Jones polynomials both for knots with high braid indices and for higher colored representations. Our data set thus comprises knots with a braid index of at most $7$.

\subsection{Jones Polynomial ($N=2$, Six-Vertex Model)}

The two-state vertex model ($N=2$) has six different configurations on the edges surrounding a vertex, or equivalently six independent Boltzmann weights at each vertex. The model is thus also called the six-vertex model. The braided $R$-matrix obtained from this model is
\begin{equation}
R = \left(
\begin{array}{cccc}
 1 & 0 & 0 & 0 \\
 0 & 0 & -q^{1/2} & 0 \\
 0 & -q^{1/2} & 1-q & 0 \\
 0 & 0 & 0 & 1 \\
\end{array}
\right).
\end{equation}

\paragraph{Two-strand braids:} This $R$-matrix allows direct computation of polynomials for knots and links represented by two-strand braids such as the unknot, trefoil, and \(T(2, 2n)\) links. For the case of the trefoil, whose braid word is \(b_1^3\), we can write the matrix for \(b_1^3\) and multiply it with $H=h \otimes h$, which for two strands takes the form
\begin{equation}
H = h \otimes h = \begin{pmatrix}
 \frac{1}{(1+q)^2} & 0 & 0 & 0 \\
 0 & \frac{q}{(1+q)^2} & 0 & 0 \\
 0 & 0 & \frac{q}{(1+q)^2} & 0 \\
 0 & 0 & 0 & \frac{q^2}{(1+q)^2} \\
\end{pmatrix} \,.
\end{equation}
Taking the trace of the resulting \(4 \times 4\) matrix gives the Markov trace \(\phi(b_1^3)\). The final answer for the trefoil is:
\begin{equation}
\alpha(b_1^3) = q + q^3 - q^4 \,.
\end{equation}

\paragraph{Three-strand braids:} For braids with three strands, we need to take tensor products of the $R$-matrix with the identity operator on the other strands. This gives us the $8\times8$ generator matrices $b_1 = R \otimes I_{2}$ and $b_2 = I_{2} \otimes R$ for the braid group $B_3$, with explicit expressions
\begin{equation}
	b_{1}=\left(
	\begin{array}{cccccccc}
	 1 & 0 & 0 & 0 & 0 & 0 & 0 & 0 \\
	 0 & 1 & 0 & 0 & 0 & 0 & 0 & 0 \\
	 0 & 0 & 0 & 0 & -\sqrt{q} & 0 & 0 & 0 \\
	 0 & 0 & 0 & 0 & 0 & -\sqrt{q} & 0 & 0 \\
	 0 & 0 & -\sqrt{q} & 0 & 1-q & 0 & 0 & 0 \\
	 0 & 0 & 0 & -\sqrt{q} & 0 & 1-q & 0 & 0 \\
	 0 & 0 & 0 & 0 & 0 & 0 & 1 & 0 \\
	 0 & 0 & 0 & 0 & 0 & 0 & 0 & 1 \\
	\end{array}
	\right)\,,\quad
	b_{2}= \left(
	\begin{array}{cccccccc}
	 1 & 0 & 0 & 0 & 0 & 0 & 0 & 0 \\
	 0 & 0 & -\sqrt{q} & 0 & 0 & 0 & 0 & 0 \\
	 0 & -\sqrt{q} & 1-q & 0 & 0 & 0 & 0 & 0 \\
	 0 & 0 & 0 & 1 & 0 & 0 & 0 & 0 \\
	 0 & 0 & 0 & 0 & 1 & 0 & 0 & 0 \\
	 0 & 0 & 0 & 0 & 0 & 0 & -\sqrt{q} & 0 \\
	 0 & 0 & 0 & 0 & 0 & -\sqrt{q} & 1-q & 0 \\
	 0 & 0 & 0 & 0 & 0 & 0 & 0 & 1 \\
	\end{array}
	\right)\,,
\end{equation}
and $H = h \otimes h \otimes h$ given by
\begin{equation}
	H = \frac1{(1+q)^3}\,\text{diag}\,(1,q,q,q^2,q,q^2,q^2,q^3)\,.
\end{equation}
The four-crossing figure-eight knot is represented by the braid word $A = b_1 b_2^{-1} b_1 b_2^{-1}\in B_3$.
It is easy to see that this leads to the fundamental Jones polynomial for the figure-eight knot,
\begin{equation} 
\alpha(4_1)=1+q^{-2}-q^{-1}-q+q^2\,.
\end{equation}

\subsection{Adjoint Jones Polynomials ($N=3$, $19$-vertex model)}

The three-state vertex model has $19$ independent Boltzmann weights at each vertex.  The $R$-matrix obtained from this model is
\begin{equation}
	R =\left(
	\begin{array}{ccccccccc}
	 1 & 0 & 0 & 0 & 0 & 0 & 0 & 0 & 0 \\
	 0 & 0 & 0 & -q & 0 & 0 & 0 & 0 & 0 \\
	 0 & 0 & 0 & 0 & 0 & 0 & q^2 & 0 & 0 \\
	 0 & -q & 0 & 1-q^2 & 0 & 0 & 0 & 0 & 0 \\
	 0 & 0 & 0 & 0 & q & 0 & -\sqrt{q}+q^{5/2} & 0 & 0 \\
	 0 & 0 & 0 & 0 & 0 & 0 & 0 & -q & 0 \\
	 0 & 0 & q^2 & 0 & -\sqrt{q}+q^{5/2} & 0 & 1-q-q^2+q^3 & 0 & 0 \\
	 0 & 0 & 0 & 0 & 0 & -q & 0 & 1-q^2 & 0 \\
	 0 & 0 & 0 & 0 & 0 & 0 & 0 & 0 & 1 \\
	\end{array}
	\right).
\end{equation}

Just as we derived the braid generators \( b_i \)'s from the six-vertex model, we can also obtain the braid generators from the $19$-vertex model. From the matrix \( h \), 
\begin{equation}
	h = \frac{1}{1+q+q^2} \left( 
	\begin{array}{ccc}
	1 & 0 & 0 \\
	0 & q & 0 \\
	0 & 0 & q^2 \\
	\end{array}
	\right) \,,
\end{equation}
we can construct the matrices $H^{(n)}=h\otimes\cdots\otimes h$ for braid words with $n$ strands. Braid generators are formed by considering tensor products such as $b_1=R\otimes I_3\otimes\cdots\otimes I_3$, etc. For example, one can calculate the following polynomials,
\begin{eqnarray}
	\alpha(3_1) &=&\ q^2 (1 + q^3 - q^5 + q^6 - q^7 - q^8 + q^9)\,,\\
	\alpha(4_1) &=&\ 3+q^{-6} - q^{-5} -q^{-4} + 2q^{-3} - q^{-2} - q^{-1}- q - q^2 +  2 q^3 - q^4 - q^5 + q^6 \,. \nonumber \\
\end{eqnarray}

\paragraph{Remark:}
The choice of framing is different here as compared to the choice in the \texttt{KnotTheory} package for Mathematica.
Let us denote a braid word, \textit{e.g.}, $A=b_ib_jb_k^{-1}$,  by a numerical list, $\tilde A=\{i,j,-k\}$.
To match the framing factor, we introduce a prefactor \( q^{s(A)-\s(A)} \), where \( s(A) \) denotes the sum of signs of exponents of generators in braid word $A$, and \( \s(A) \) denotes the sum of the numerical braid word representation.
As an example, for the braid word \( A= b_1 b_2^{-1} b_1 b_2^{-1} \), or equivalently \(\tilde A= \{1, -2, 1, -2\} \), these are
\begin{equation}
s(A) = 1 - 1 + 1 - 1 = 0 \,, \quad \s(A) = 1 - 2 + 1 - 2 = -2 \,.
\end{equation}
In this example, the prefactor becomes \( q^{0 - (-2)} = q^2 \). Multiplying the vertex model answer with this prefactor, we get a polynomial that matches exactly with that obtained from the \texttt{KnotTheory} package.

\bibliographystyle{JHEP}
\bibliography{ref}

\end{document}

%% file: custom.tex
\newcommand{\Rmnum}[1]{\expandafter\@slowromancap\romannumeral #1@}

\newcommand{\lB}{\left [}
\newcommand{\rB}{\right ]}
\newcommand{\lb}{\left (}
\newcommand{\rb}{\right )}
\renewcommand{\a}{\alpha}	%%% Redefinition
\newcommand{\f}{\phi}

\newcommand{\s}{\sigma}

\renewcommand{\t}{\tau}		%%% Redefinition
\newcommand{\I}{\infty}

%% file: main.bbl
\providecommand{\href}[2]{#2}\begingroup\raggedright\begin{thebibliography}{10}

\bibitem{Witten:1988hf}
E.~Witten, \emph{{Quantum Field Theory and the Jones Polynomial}},
  \href{http://dx.doi.org/10.1007/BF01217730}{\emph{Commun. Math. Phys.} {\bf
  121} (1989) 351--399}.

\bibitem{jones1985polynomial}
V.~F. Jones, \emph{A polynomial invariant for knots via von neumann algebras},
  {\emph{Bulletin of the American Mathematical Society} {\bf 12} (1985)
  103--111}.

\bibitem{jones1987hecke}
V.~Jones, \emph{Hecke algebra representations of braid groups and link
  polynomials}, {\emph{Annals of Mathematics} (1987) 335--388}.

\bibitem{kauffman1987state}
L.~H. Kauffman, \emph{State models and the jones polynomial}, {\emph{Topology}
  {\bf 26} (1987) 395--407}.

\bibitem{khovanov2000categorification}
M.~Khovanov, \emph{A categorification of the jones polynomial}, {\emph{Duke
  Math. J.} {\bf 104} (2000) 359--426}.

\bibitem{bar2002khovanov}
D.~Bar-Natan, \emph{On khovanov’s categorification of the jones polynomial},
  {\emph{Algebraic \& Geometric Topology} {\bf 2} (2002) 337--370}.

\bibitem{thurston_geometry_topology_three_manifolds}
W.~P. Thurston, \emph{The Geometry and Topology of Three-Manifolds}.
\newblock Princeton University Press, 1979.

\bibitem{belousov2019hyperbolicknotsgeneric}
Y.~Belousov and A.~Malyutin, \emph{Hyperbolic knots are not generic},  2019.

\bibitem{kashaev1995link}
R.~M. Kashaev, \emph{A link invariant from quantum dilogarithm}, {\emph{Modern
  Physics Letters A} {\bf 10} (1995) 1409--1418}.

\bibitem{murakami2001colored}
H.~Murakami and J.~Murakami, \emph{The colored jones polynomials and the
  simplicial volume of a knot}, {\emph{Acta Mathematica} {\bf 186} (2001)
  85--104}.

\bibitem{Gukov:2003na}
S.~Gukov, \emph{{Three-dimensional quantum gravity, Chern-Simons theory, and
  the A polynomial}},
  \href{http://dx.doi.org/10.1007/s00220-005-1312-y}{\emph{Commun. Math. Phys.}
  {\bf 255} (2005) 577--627}, [\href{https://arxiv.org/abs/hep-th/0306165}{{\tt
  hep-th/0306165}}].

\bibitem{witten2011analytic}
E.~Witten, \emph{Analytic continuation of chern-simons theory}, {\emph{AMS/IP
  Stud. Adv. Math} {\bf 50} (2011) 347}.

\bibitem{Craven:2020bdz}
J.~Craven, V.~Jejjala and A.~Kar, \emph{{Disentangling a deep learned volume
  formula}}, \href{http://dx.doi.org/10.1007/JHEP06(2021)040}{\emph{JHEP} {\bf
  06} (2021) 040}, [\href{https://arxiv.org/abs/2012.03955}{{\tt 2012.03955}}].

\bibitem{Garoufalidis_2011}
S.~Garoufalidis and T.~Lê, \emph{Asymptotics of the colored jones function of
  a knot}, \href{http://dx.doi.org/10.2140/gt.2011.15.2135}{\emph{Geometry \&
  Topology} {\bf 15} (Oct., 2011) 2135–2180}.

\bibitem{Habiro_2007}
K.~Habiro, \emph{A unified witten–reshetikhin–turaev invariant for integral
  homology spheres},
  \href{http://dx.doi.org/10.1007/s00222-007-0071-0}{\emph{Inventiones
  mathematicae} {\bf 171} (Sept., 2007) 1–81}.

\bibitem{KnotAtlas}
``The knot atlas.'' URL: \url{katlas.org/}.

\bibitem{linkinfo}
C.~Livingston and A.~H. Moore, ``Linkinfo: Table of link invariants.'' URL:
  \url{linkinfo.math.indiana.edu}.

\bibitem{SnapPy}
M.~Culler, N.~M. Dunfield, M.~Goerner and J.~R. Weeks, ``Snap{P}y, a computer
  program for studying the geometry and topology of $3$-manifolds.'' Available
  at \url{http://snappy.computop.org},.

\bibitem{Jejjala:2019kio}
V.~Jejjala, A.~Kar and O.~Parrikar, \emph{{Deep Learning the Hyperbolic Volume
  of a Knot}},
  \href{http://dx.doi.org/10.1016/j.physletb.2019.135033}{\emph{Phys. Lett. B}
  {\bf 799} (2019) 135033}, [\href{https://arxiv.org/abs/1902.05547}{{\tt
  1902.05547}}].

\bibitem{Le:2000}
T.~T.~Q. L{\^e}, \emph{Quantum invariants of 3-manifolds: Integrality,
  splitting, and perturbative expansion},
  \href{http://dx.doi.org/10.1016/S0166-8641(02)00056-1}{\emph{Topology and its
  Applications} {\bf 127} (2000) 125--152},
  [\href{https://arxiv.org/abs/math/0004099}{{\tt math/0004099}}].

\bibitem{Habiro:2000}
K.~Habiro, \emph{On the colored jones polynomials of some simple links},
  {\emph{Surikaisekikenkyusho Kokyuroku} {\bf 1172} (2000) 34}.

\bibitem{Garoufalidis:2004}
S.~Garoufalidis and Y.~Lan, \emph{Experimental evidence for the volume
  conjecture for the simplest hyperbolic non-2–bridge knot},
  \href{http://dx.doi.org/10.2140/agt.2005.5.379}{\emph{Algebraic \& Geometric
  Topology} {\bf 5} (2004) 379--403},
  [\href{https://arxiv.org/abs/math/0412331}{{\tt math/0412331}}].

\bibitem{hoste1998first}
J.~Hoste, M.~Thistlethwaite and J.~Weeks, \emph{The first 1,701,936 knots},
  {\emph{Math. Intelligencer} {\bf 20} (1998) 33--48}.

\bibitem{Akutsu:1987dz}
Y.~Akutsu and M.~Wadati, \emph{{Knot invariants and the critical statistical
  systems}}, \href{http://dx.doi.org/10.1143/JPSJ.56.839}{\emph{J. Phys. Soc.
  Jap.} {\bf 56} (1987) 839--842}.

\bibitem{Ramadevi2017}
P.~Ramadevi, \emph{Vertex Models and Knot invariants}.
\newblock Springer Singapore, Singapore, 2017,
  \href{http://dx.doi.org/10.1007/978-981-10-6841-6\_14}{10.1007/978-981-10-6841-6\_14}.

\bibitem{hughes2016neural}
M.~C. Hughes, \emph{A neural network approach to predicting and computing knot
  invariants}, {\emph{Journal of Knot Theory and Its Ramifications} {\bf 29}
  (2020) 2050005}, [\href{https://arxiv.org/abs/1610.05744}{{\tt 1610.05744}}].

\bibitem{levitt2022big}
J.~S. Levitt, M.~Hajij and R.~Sazdanovic, \emph{Big data approaches to knot
  theory: understanding the structure of the jones polynomial}, {\emph{Journal
  of Knot Theory and Its Ramifications} {\bf 31} (2022) 2250095},
  [\href{https://arxiv.org/abs/1912.10086}{{\tt 1912.10086}}].

\bibitem{Gukov:2020qaj}
S.~Gukov, J.~Halverson, F.~Ruehle and P.~Sulkowski, \emph{{Learning to
  Unknot}}, \href{http://dx.doi.org/10.1088/2632-2153/abe91f}{\emph{Mach.
  Learn. Sci. Tech.} {\bf 2} (2021) 025035},
  [\href{https://arxiv.org/abs/2010.16263}{{\tt 2010.16263}}].

\bibitem{dlotko2023mappertypealgorithmscomplexdata}
P.~Dlotko, D.~Gurnari and R.~Sazdanovic, \emph{Mapper-type algorithms for
  complex data and relations},  \href{https://arxiv.org/abs/2109.00831}{{\tt
  2109.00831}}.

\bibitem{davies2021advancing}
A.~Davies, P.~Veli{\v{c}}kovi{\'c}, L.~Buesing, S.~Blackwell, D.~Zheng,
  N.~Toma{\v{s}}ev, R.~Tanburn et~al., \emph{Advancing mathematics by guiding
  human intuition with ai}, {\emph{Nature} {\bf 600} (2021) 70--74}.

\bibitem{Craven:2021ckk}
J.~Craven, M.~Hughes, V.~Jejjala and A.~Kar, \emph{{Learning knot invariants
  across dimensions}},
  \href{http://dx.doi.org/10.21468/SciPostPhys.14.2.021}{\emph{SciPost Phys.}
  {\bf 14} (2023) 021}, [\href{https://arxiv.org/abs/2112.00016}{{\tt
  2112.00016}}].

\bibitem{Craven:2022cxe}
J.~Craven, M.~Hughes, V.~Jejjala and A.~Kar, \emph{{Illuminating new and known
  relations between knot invariants}},
  \href{http://dx.doi.org/10.1088/2632-2153/ad95d9}{\emph{Mach. Learn. Sci.
  Tech.} {\bf 5} (2024) 045061}, [\href{https://arxiv.org/abs/2211.01404}{{\tt
  2211.01404}}].

\bibitem{Gukov:2023kvx}
S.~Gukov, J.~Halverson, C.~Manolescu and F.~Ruehle, \emph{{Searching for
  ribbons with machine learning}},
  \href{https://arxiv.org/abs/2304.09304}{{\tt 2304.09304}}.

\bibitem{Gukov:2024buj}
S.~Gukov, J.~Halverson and F.~Ruehle, \emph{{Rigor with machine learning from
  field theory to the Poincar\'e conjecture}},
  \href{http://dx.doi.org/10.1038/s42254-024-00709-0}{\emph{Nature Rev. Phys.}
  {\bf 6} (2024) 310--319}, [\href{https://arxiv.org/abs/2402.13321}{{\tt
  2402.13321}}].

\bibitem{Gukov:2024opc}
S.~Gukov and R.-K. Seong, \emph{{Machine learning BPS spectra and the gap
  conjecture}},
  \href{http://dx.doi.org/10.1103/PhysRevD.110.046016}{\emph{Phys. Rev. D} {\bf
  110} (2024) 046016}, [\href{https://arxiv.org/abs/2405.09993}{{\tt
  2405.09993}}].

\bibitem{github}
M.~Hughes, V.~Jejjala, P.~Ramadevi, P.~Roy and V.~K. Singh,
  ``{Colored}{Jones}{ML}.''
  \url{https://github.com/roypratik92/ColoredJonesML}, Feb., 2025.

\bibitem{data}
M.~Hughes, V.~Jejjala, P.~Ramadevi, P.~Roy and V.~K. Singh, ``Data for 2- and
  3-dimensional colored jones polynomials for hyperbolic knots.''
  \url{https://doi.org/10.5281/zenodo.1490093}, Feb., 2025.
\newblock 10.5281/zenodo.14900936.

\bibitem{yokota2000volume}
Y.~Yokota, \emph{On the volume conjecture for hyperbolic knots}, {\emph{arXiv
  preprint math/0009165} (2000) }.

\bibitem{Ishikawa2019}
K.~Ishikawa, \emph{Hoste’s conjecture for the 2-bridge knots},
  {\emph{Proceedings of the American Mathematical Society} (2019) }.

\bibitem{Stoimenow2011}
A.~{Stoimenow}, \emph{{Diagram genus, generators and applications}},
  \href{http://dx.doi.org/10.48550/arXiv.1101.3390}{\emph{arXiv e-prints}
  (Jan., 2011) arXiv:1101.3390}, [\href{https://arxiv.org/abs/1101.3390}{{\tt
  1101.3390}}].

\bibitem{alsukaiti2024alexander}
M.~E. AlSukaiti and N.~Chbili, \emph{Alexander and jones polynomials of weaving
  3-braid links and whitney rank polynomials of lucas lattice}, {\emph{Heliyon}
  {\bf 10} (2024) }.

\bibitem{HM2013}
M.~Hirasawa and K.~Murasugi, \emph{Various stabilities of the alexander
  polynomials of knots and links}, {\emph{arXiv preprint arXiv:1307.1578}
  (2013) }.

\bibitem{HIRASAWA201948}
M.~Hirasawa, K.~Ishikawa and M.~Suzuki, \emph{Alternating knots with alexander
  polynomials having unexpected zeros},
  \href{http://dx.doi.org/https://doi.org/10.1016/j.topol.2018.11.027}{\emph{Topology
  and its Applications} {\bf 253} (2019) 48--56}.

\bibitem{Jin2010}
F.~Zhang, F.~Dong, E.~G. Tay et~al., \emph{Zeros of the jones polynomial are
  dense in the complex plane}, {\emph{the electronic journal of combinatorics}
  (2010) R94--R94}.

\bibitem{AI2005}
A.~Champanerkar and I.~Kofman, \emph{On the mahler measure of jones polynomials
  under twisting}, {\emph{Algebraic \& Geometric Topology} {\bf 5} (2005)
  1--22}.

\bibitem{XF2010}
X.~{Jin} and F.~{Zhang}, \emph{{Zeros of the Jones Polynomial for Multiple
  Crossing-Twisted Links}},
  \href{http://dx.doi.org/10.1007/s10955-010-0027-4}{\emph{Journal of
  Statistical Physics} {\bf 140} (Sept., 2010) 1054--1064}.

\bibitem{Andersen:2004}
J.~E. Andersen, N.~Askitas, D.~Bar-Natan, S.~Baseilhac, R.~Benedetti,
  S.~Bigelow, M.~Boileau et~al., \emph{Problems on invariants of knots and
  3-manifolds}, {\emph{Geometry \& Topology Monographs} {\bf 4} (2002)
  377--572}, [\href{https://arxiv.org/abs/math/0406190}{{\tt math/0406190}}].

\bibitem{williamson2024deep}
G.~Williamson, \emph{Is deep learning a useful tool for the pure
  mathematician?}, {\emph{Bulletin of the American Mathematical Society} {\bf
  61} (2024) 271--286}, [\href{https://arxiv.org/abs/2304.12602}{{\tt
  2304.12602}}].

\bibitem{Murakami:2022}
H.~Murakami, \emph{The colored jones polynomial of the figure-eight knot and a
  quantum modularity},
  \href{http://dx.doi.org/10.4153/S0008414X23000172}{\emph{Canadian Journal of
  Mathematics} {\bf 76} (2022) 519 -- 554},
  [\href{https://arxiv.org/abs/2209.07751}{{\tt 2209.07751}}].

\bibitem{KirResh}
A.~Kirillov and N.~Reshetikhin, \emph{Representations of the Algebra
  $U_q(sl_2)$, q-Orthogonal Polynomials and Invariants of Links}.
\newblock World Scientific, 1990.

\bibitem{RT1}
E.~Guadagnini, M.~Martellini and M.~Mintchev, \emph{Quantum groups},  in
  \emph{New Developments in the Theory of Knots}, pp.~307--317, World
  Scientific, 1989.

\bibitem{RT2}
N.~Reshetikhin and V.~Turaev, \emph{Ribbon graphs and their invariants derived
  from quantum groups}, {\emph{Comm. Math. Phys.} {\bf 127} (1990) 1--26}.

\bibitem{Itoyama:2012qt}
H.~Itoyama, A.~Mironov, A.~Morozov and A.~Morozov, \emph{Character expansion
  for homfly polynomials. iii. all 3-strand braids in the first symmetric
  representation},
  \href{http://dx.doi.org/10.1142/S0217751X12500996}{\emph{Int. J. Mod. Phys.
  A} {\bf 27} (2012) 1250099}, [\href{https://arxiv.org/abs/arXiv:1204.4785
  [hep-th]}{{\tt arXiv:1204.4785 [hep-th]}}].

\bibitem{Kaul:1992rs}
R.~K. Kaul and T.~R. Govindarajan, \emph{{Three-dimensional Chern-Simons theory
  as a theory of knots and links. 2. Multicolored links}},
  \href{http://dx.doi.org/10.1016/0550-3213(93)90251-J}{\emph{Nucl. Phys.} {\bf
  B393} (1993) 392--412}.

\bibitem{Devi:1993ue}
P.~R. Devi, T.~R. Govindarajan and R.~K. Kaul, \emph{Three dimensional
  chern-simons theory as a theory of knots and links iii : Compact semi-simple
  group},
  \href{http://dx.doi.org/10.1016/0550-3213(93)90652-6}{\emph{Nucl.Phys.} {\bf
  B402} (1993) 548--566}, [\href{https://arxiv.org/abs/hep-th/9212110}{{\tt
  hep-th/9212110}}].

\bibitem{RGK}
P.~Ramadevi, T.~R. Govindarajan and R.~K. Kaul, \emph{{Knot invariants from
  rational conformal field theories}},
  \href{http://dx.doi.org/10.1016/0550-3213(94)00102-2}{\emph{Nucl. Phys.} {\bf
  B422} (1994) 291--306}, [\href{https://arxiv.org/abs/hep-th/9312215}{{\tt
  hep-th/9312215}}].

\bibitem{DMMMRSS:2017}
S.~Dhara, A.~Mironov, A.~Morozov, A.~Morozov, P.~Ramadevi, V.~K. Singh and
  A.~Sleptsov, \emph{{Eigenvalue hypothesis for multistrand braids}},
  \href{http://dx.doi.org/10.1103/PhysRevD.97.126015}{\emph{Phys. Rev. D} {\bf
  97} (2018) 126015}, [\href{https://arxiv.org/abs/1711.10952}{{\tt
  1711.10952}}].

\bibitem{DMMMRSS:2018}
S.~Dhara, A.~Mironov, A.~Morozov, A.~Morozov, P.~Ramadevi, V.~K. Singh and
  A.~Sleptsov, \emph{{Multi-Colored Links From 3-strand Braids Carrying
  Arbitrary Symmetric Representations}},
  \href{http://dx.doi.org/10.1007/s00023-019-00841-z}{\emph{Annales Henri
  Poincare} {\bf 20} (2019) 4033--4054},
  [\href{https://arxiv.org/abs/1805.03916}{{\tt 1805.03916}}].

\bibitem{Wadati:1989ud}
M.~Wadati, T.~Deguchi and Y.~Akutsu, \emph{{Exactly Solvable Models and Knot
  Theory}}, \href{http://dx.doi.org/10.1016/0370-1573(89)90123-3}{\emph{Phys.
  Rept.} {\bf 180} (1989) 247}.

\bibitem{Baxter:1982zz}
R.~J. Baxter, \emph{{Exactly solved models in statistical mechanics}}.
\newblock 1982,
  \href{http://dx.doi.org/10.1142/9789814415255\_0002}{10.1142/9789814415255\_0002}.

\bibitem{Akutsu:1987qs}
Y.~Akutsu and M.~Wadati, \emph{{Exactly Solvable Models and New Link
  Polynomials. 1. $N$ State Vertex Models}},
  \href{http://dx.doi.org/10.1143/JPSJ.56.3039}{\emph{J. Phys. Soc. Jap.} {\bf
  56} (1987) 3039--3051}.

\end{thebibliography}\endgroup
